\documentclass[12pt]{article}
\usepackage{amssymb}
\usepackage{amsmath}
\usepackage{graphicx,color}
\usepackage{amsfonts}
\usepackage[margin=1in]{geometry}
\usepackage{url}
\usepackage{natbib} 
\usepackage{xcolor}
\usepackage[colorlinks,citecolor=blue!65!black,linkcolor=black,urlcolor=blue]{hyperref}
\usepackage{subcaption}
\usepackage{xr}
\usepackage{sectsty}
\usepackage{tikz}
\usepackage{pgfplots}
\sectionfont{\centering} 
\subsectionfont{\centering} 
\usepackage{setspace}
 
\usepackage{caption}
\usepackage{subcaption}

\usepackage{comment}
\providecommand{\keywords}[1]{{\hspace{1em}\textit{Keywords:}} #1}
\usepackage{abstract}

\newcommand{\RNum}[1]{\uppercase\expandafter{\romannumeral #1\relax}}
\newcommand{\md}{\mathrm{d}}

\newtheorem{theorem}{Theorem}

\newtheorem{condition}{Condition}

\newtheorem{lemma}{Lemma}

\newtheorem{proposition}{Proposition}
\newtheorem{remark}{Remark}

\newcommand{\indicator}{\mathrm{I}}

\title{Extended Asymptotic Identifiability of Nonparametric Item Response Models}
\author{Yinqiu He}
\date{Department of Statistics\\ University of Wisconsin-Madison}
\begin{document}
\maketitle

\begin{abstract}
Nonparametric item response    models provide a flexible   framework in psychological and educational measurements.  
\cite{douglas2001asymptotic} established  asymptotic identifiability for a class of models with nonparametric response functions  for long assessments.  
Nevertheless,  the model class examined in \cite{douglas2001asymptotic} excludes several popular parametric item response models. 
This limitation can hinder the applications in which  nonparametric and parametric models   are compared, such as evaluating  model  goodness-of-fit.   
To address this issue, 
We consider an extended nonparametric model class   that encompasses most parametric models and establish asymptotic identifiability.  
The results bridge the parametric and nonparametric item response models and provide  a solid theoretical foundation for the applications of nonparametric item response models for assessments with many items.
\end{abstract}


\medskip
\section{Introduction}
Item response theory (IRT) models play a crucial role in 
   psychological measurements, educational testing, and political science \citep{van2018handbook}. 
The parametric IRT models include a broad class of models that assume parametric forms of item response  functions (IRFs),  such as  the normal-ogive models and the three-parameter logistic models. 
%In a typical set-up of educational testing, dichotomous responses to $n$ items $(Y_1,\ldots, Y_n)$ are observed and 
%psychometricians are interested in studying the 
Nevertheless, it  has been widely   recognized that  IRFs  cannot always be modeled well with parametric families. This has spurred extensive research into the theory and applications of nonparametric IRT  models; see  reviews in  \cite{sijtsma1998methodology},  \cite{sijtsma2002introduction}, and \cite{chen2021item}. 
Nonparametric IRT models have been popularly used 
%play an important role 
in assessing the goodness-of-fit of parametric IRT models and in providing robust measurements against model misfitting.

In the pursuit of  nonparametric modelling of IRT models,   
one pioneer is the  
Mokken scale analysis  
\citep{mokken1982nonparametric}.  
This analysis is designed to partition a set of items into Mokken scales and can also be employed to  check   assumptions of nonparametric IRT models \citep{van2007mokken}. 
Mokken's work was further developed by \cite{sijtsma1987reliability},  
\cite{sijtsma1998methodology} and \cite{sijtsma2017tutorial},  among    others. 
Independent of this research line,
another strand of research has  focused on   modeling   nonparametric IRFs. 
These studies relax  the assumptions of parametric   IRFs via nonparametric functions such as  splines or polynomials 
\citep{winsberg1984fitting,ramsay1989binomial,ramsay1991maximum,ramsay1991kernel,douglas1997joint,johnson2007modeling,peress2012identification,falk2016maximum}. 
%semiparametric models in \cite{peress2012identification}. 

%\cite{junker1997characterization}
%efforts in all kinds of identifiability results 

%semiparametric \cite{cressie1983characterizing} \cite{lindsay1991semiparametric}  \cite{ramsay1991maximum} \cite{peress2012identification}
%other: \cite{johnson2006nonparametric} other discussions: future work latent variable dimension (univariate) \cite{peress2012identification}

When the true IRFs are assumed to belong to a very general function space, 
there may exist different sets of distinct IRFs that yield identical distributions of manifest variables. 
Understanding the identifiability of nonparametric IRT models is critical for relating the obtained estimates to the underlying true models.  
To address this issue, 
\cite{douglas2001asymptotic} established the identifiability of nonparametric IRT models 
in an asymptotic sense with the number of items $n$ going to infinity.
The theoretical results provide foundations for various applications of nonparametric IRTs, including assessing the   parametric model fit \citep{douglas2001nonparametric,lee2009use}. 
% motivates research on  the goodness-of-fit \cite{douglas2001nonparametric} and \cite{lee2009use}. 
%\paragraph{review} of \cite{douglas2001asymptotic} 

Nevertheless, as pointed out in  \cite{douglas2001asymptotic}, their identifiability result relies on restrictive assumptions about the model class of IRFs, which are not met by some popular parametric item response models such as the normal ogive model. 
But in applications such as assessing the parametric model fit, it is often desired to consider a class of  nonparametric IRTs that can encompass  the widely-used parametric IRTs.
For example,  \cite{lee2009use} proposed to evaluate the fit of a two-parameter logistic  (2PL) model by comparing  the estimated 2PL with another estimated nonparametric IRT. 
However, as the existing  identifiability results in \cite{douglas2001asymptotic} do  not include 2PL in the model class, 
there may exist a nonparametric IRF that differs from the 2PL IRF but  yields the same manifest distribution. 
Consequently, the discrepancy between the nonparametric and the 2PL IRFs might be a trivial result from  non-identifiability rather than model misfit.
This could lead to unreliability in goodness-of-fit measurements based on nonparametric IRTs. 
%, so that the goodness-of-fit measurements based on nonparametric IRTs  may be unreliable.    
 
%Without a   clear understanding of the identifiability results covering the 

To lay a solid theoretical foundation for  related applications, 
it is imperative to establish identifiability for a model class that extends that in  \cite{douglas2001asymptotic}. 
%In brief, establishing the   identifiability for a model class that extends the scope of \cite{douglas2001asymptotic} would provide solid theoretical foundations in the related applications. 
In this paper, 
we relax the assumptions made in \cite{douglas2001asymptotic} so that the model class can encompass a broader range of parametric item response models, such as the normal-ogive model and the four-parameter logistic model. 
%This may limit the use of  the nonparametric IRTs. 

Following \cite{douglas2001asymptotic}, we focus on monotone, unidimensional, and locally independent item response models.  
Specifically, 
let $\boldsymbol{Y}_n=(Y_1,\ldots, Y_n)$ represent $n$ observed dichotomous variables, and 
denote the $i$-th IRF as $P_{i}(\theta)=P(Y_i=1\mid \Theta=\theta)\in [0,1]$, which is unidimensional and strictly increasing for $i=1,\ldots, n$. Here, $P_i(\theta)$   represents the probability that a person with a given ability level $\theta$ will answer the $i$-th item  correctly.  
Then for the manifest distribution of $\boldsymbol{Y}_n$, 
%a monotone, unidimensional, and locally independent item response model consists of a collection of strictly increasing probability functions 
a locally independent IRT model consists of a collection of 
functions $\{P_1, P_2,\ldots, P_n\}$ and a probability density function $f$ of $\theta$ that satisfy 
%A monotone, unidimensional, and locally independent item response model for the manifest distribution of $n$ dichotomous variables, $(Y_1,\ldots, Y_n)$, consists of a  collection of strictly increasing probability distribution functions $\{P_1, P_2,\ldots, P_n\}$ taking values between 0 and 1, and a probability density function $f$ of $\theta$  satisfying,
\begin{align}\label{eq:integralform}
	 P\left[Y_{i_1}=1, Y_{i_2}=1, \ldots, Y_{i_k}=1\right]=\int \prod_{j=1}^k P_{i_j}(\theta) f(\theta) d \theta 
\end{align}
for any nonempty subsets $\left\{i_1, i_2, \ldots, i_k\right\}$ of the integers $\{1,2, \ldots, n\}$. 
Given a fixed probability density function  $f(\theta)$, we say that the IRFs are identifiable, if for any other collection of IRFs $\left\{P_1^*, P_2^*, \ldots, P_n^*\right\}$ satisfying \eqref{eq:integralform}, we have $P_i=P_i^*$ for all $i \in\{1,2, \ldots, n\}$.

%\cite{douglas2001asymptotic} established the asymptotic identifiability with $n\to \infty$ assuming that the derivatives of the IRFs are uniformly bounded away from 0 and $\infty$, and the  limits of IRFs at the ends points are 0 or 1. 

It is worth noting that the analyses in this paper and  \cite{douglas2001asymptotic} restrict $f$ to be a fixed density function. 
Allowing the transformation of $f$ can 
introduce non-identifiability issues that may not contribute meaningfully to the analysis. 
 To illustrate,  consider a transformed latent trait $\lambda=G(\theta)$ for a monotone function $G(\cdot)$, and  $G_{inv}(\cdot)$, the inverse function  of $G(\cdot)$, satisfies $\theta = G_{inv}(\lambda)$.  
Then the manifest distribution  \eqref{eq:integralform} is equivalent to
\begin{align*}
 P\left[Y_{i_1}=1, Y_{i_2}=1, \ldots, Y_{i_k}=1\right]=\int \prod_{j=1}^k P_{i_j}[G^{-1}(\lambda)] f[G^{-1}(\lambda)] G'_{inv}(\lambda) d \lambda,   	
\end{align*}
which gives rise to another IRT model with a different set of IRFs $\{P_{1}[G^{-1}(\lambda)], \ldots, P_{n}[G^{-1}(\lambda)] \}$ coupled with the latent trait density $f[G^{-1}(\lambda)] G'_{inv}(\lambda)$. 
As the choice of $f$
  can be arbitrary, we consider a fixed   $f$  for theoretical convenience.  In Section  \ref{sec:example}, we will show that our results encompass the commonly used normal distribution of the latent trait with a proper transformation. 
 
%Some other studies consider fixing the IRFs and varying the distribution of the latent trait \citep{cressie1983characterizing,lindsay1991semiparametric,haberman2005identifiability} 
%\cite{lindsay1991semiparametric},   and 
%\cite{haberman2005identifiability}

%    models where the IRFs are fixed but the   distribution of the latent trait is non-parametric  
%\cite{cressie1983characterizing} \cite{lindsay1991semiparametric},   and 
%\cite{haberman2005identifiability}. 

%We relax the assumptions in \cite{douglas2001asymptotic} so that the model class can cover a wide range of parametric item response models, including normal-ogive model and three parameter model. 
%Being able to achieve this,
%provide solid theoretical foundation that 
% pave the way for the use of assessing the goodness-of-fit \cite{douglas2001nonparametric,lee2009use}. 

%splines: \cite{johnson2007modeling}, \cite{winsberg1984fitting}%\cite{johnson2006nonparametric}
%only the IRFs are considered non-parametric
%\cite{sijtsma2002introduction}

In Section \ref{sec:setup}, we introduce an analysis framework of triangular sequences  similarly to \cite{douglas2001asymptotic}.  
Under this framework, we introduce conditions on the nonparametric item response models and present examples to demonstrate that the new conditions can significantly relax assumptions in \cite{douglas2001asymptotic}. 
Section \ref{sec:results} presents the  asymptotic identifiability results under relaxed assumptions and proofs. 
Section \ref{sec:disc} discusses the  practical implications of the results. 
Then appendix provides technical lemmas and the proofs of  lemmas and propositions. 

\section{Set-up}\label{sec:setup}
%In this paper, 
%we consider the triangular sequences of item response vectors and item characteristics curves similarly to \cite{douglas2001asymptotic}. 
We consider the triangular sequence of item response variables that can be expressed as 
\begin{align*}
	 \begin{aligned} & \boldsymbol{Y}_k=\left(Y_{k, 1}, Y_{k, 2}, \ldots, Y_{k, k}\right) \\ & \boldsymbol{Y}_{k+1}=\left(Y_{k+1,1}, Y_{k+1,2}, \ldots, Y_{k+1, k}, Y_{k+1, k+1}\right) \\ & \boldsymbol{Y}_{k+2}=\left(Y_{k+2,1}, Y_{k+2}, \ldots, 2, Y_{k+2, k+1}, Y_{k+2, k+2}\right)\\
	 & \ldots \ldots
	\end{aligned}
\end{align*}
Item response vectors in the sequence are not required to share any items with one another. 
That is to say, 
the vectors in the sequence are allowed to be disjoint, or they may overlap to any extent. 
Moreover, 
we let $\mathcal{F}_k$ denote the probability distribution of $\boldsymbol{Y}_k$.
Then  $\{\mathcal{F}_k , \mathcal{F}_{k+1}, \ldots \}$  form a triangle sequence of  distributions of item response vectors. 
In addition, we  define  a triangular sequence of IRFs as 
\begin{align*}
\begin{aligned} & \mathcal{P}_k=\left\{P_{k, 1}, P_{k, 2}, \ldots, P_{k, k}\right\} \\ & \mathcal{P}_{k+1}=\left\{P_{k+1,1}, P_{k+1,2}, \ldots, P_{k+1, k}, P_{k+1, k+1}\right\} \\ & \mathcal{P}_{k+2}=\left\{P_{k+2,1}, P_{k+2,2}, \ldots, P_{k+2, k+1}, P_{k+2, k+2}\right\}\\
&\ldots \ldots \end{aligned}	
\end{align*}
Let $f(\theta)$ be a fixed probability density function of the latent trait $\Theta$. 
We say that  $\{\mathcal{P}_k , \mathcal{P}_{k+1}, \ldots \}$, coupled with $f$, is a model for the  sequence of manifest distributions $\{\mathcal{F}_k , \mathcal{F}_{k+1}, \ldots \}$, if for each $k$, all equations between  the  manifest distributions and the integrals, specified as in \eqref{eq:integralform}, are satisfied. 
We establish the asymptotic  identifiability in the sense  that if the two   sequences of models,  $\{\mathcal{P}_k , \mathcal{P}_{k+1}, \ldots \}$ and  $\{\mathcal{P}_k^* , \mathcal{P}_{k+1}^*, \ldots \}$, are for the same sequence of manifest distributions,
their pointwise difference converges to as the number of items $n$ increase to infinity. 

%We next examine the differences between two sequences of models,  $\{\mathcal{P}_k , \mathcal{P}_{k+1}, \ldots \}$ and  $\{\mathcal{P}_k^* , \mathcal{P}_{k+1}^*, \ldots \}$,  for the same sequence of manifest distributions.

We next state conditions that specify the class of item response models to consider. 
Then we provide examples showing that the specified model class can include a wide range of popular item response models and greatly extends the model class in \cite{douglas2001asymptotic}.

%\subsection*{Conditions}
%The   theoretical results on asymptotic identifiability are presented afterwards. 

\subsection{Conditions}

%We next present the conditions and discuss their 

\begin{condition}\label{cond:unilocal}
\textit{Unidimensionality and local independence:} The latent variable $\Theta $ is a scalar valued random variable and item responses are mutually independent conditioning on $\Theta$.  
%\textit{Local independence:}. 	

\end{condition}
 
 \begin{condition}\label{cond:uniformdist}
$\Theta$ follows $U(0,1)$, a uniform distribution on the interval $(0,1)$. 
\end{condition}

Given the unidimensionality in Condition \ref{cond:unilocal}, Condition \ref{cond:uniformdist} can be viewed as the choice of a specific  parameterization for $\Theta$.   
In particular, when $\tilde{\Theta}$  is a  random variable with a continuous cumulative distribution function $F: \mathbb{R}\to (0,1)$,  
we define $\Theta=F(\tilde{\Theta})$, resulting in $\Theta\sim U(0,1)$; see   Proposition 3.1 in \cite{embrechts2013note} for details. 
This transformation suggests that any  latent trait following a continuous distribution can be equivalently transformed to $U(0,1)$.

   \begin{condition}\label{cond:derivbound}
 Given each pair $(n,i)$, the first-order derivative $P'_{n,i}(\theta)$ exists and is continuous in the open interval $(0,1)$. 
For $0<\alpha<\beta<1$, 
there exist constants $m_{\alpha\beta}$ and $M_{\alpha\beta}$ that do not depend on $(n,i)$ such that for $\theta\in [\alpha,\beta]$,
$
	0<m_{\alpha\beta}<P'_{n,i}(\theta) < M_{\alpha\beta} < \infty. 
$ 
%Given $0<\delta<1$, there exist positive numbers $m_{\delta}$ and $M_{\delta}$ such that  for $\theta\in (\delta, 1-\delta)$ and any $n$ and $i$,  $m_{\delta}<P^{\prime}_{n,i}(\theta)<M_{\delta}$.
\end{condition}

Condition \ref{cond:derivbound} requires  that the   derivatives  $P'_{n,i}(\theta)$'s are uniformly bounded from below and above on a compact interval $[\alpha,\beta]\subseteq (0,1)$.   
This is a notable relaxation of Assumption 4 in \cite{douglas2001asymptotic}, which requires that $P'_{n,i}(\theta)$'s are uniformly bounded over the entire interval $(0,1)$. 
With $m_{\alpha\beta}>0$ in Condition \ref{cond:derivbound}, we ensure  that $P_{n,i}(\theta)$'s are   strictly  increasing with respect to $\theta$, which is a commonly accepted  assumption in the literature. 
Nevertheless, we point out that all the analyses can be readily extended to cases where   $P_{n,i}(\theta)$'s are   strictly  decreasing with respect to $\theta$ by consider the transformation $\tilde{\Theta}=1-\Theta \sim U(0,1)$ and $\tilde{P}_{n,i}(\theta)=P_{n,i}(1-\theta)$ so that $\tilde{P}_{n,i}'(\theta)=-P_{n,i}'(1-\theta)$. 

%Moreover, in some scenarios 
%\citep{sijtsma1998methodology},
%it could be desired to allow some individual IRFs to be locally decreasing. 
%Condition \ref{cond:derivbound} can be replaced and the proof would still hold. 
%\begin{remark}
%	In some scenarios \cite{sijtsma1998methodology}, it may be desired to relax the monotonicity assumption of individual IRTs.
%	Condition \ref{cond:derivbound} can be relaxed to mean IRTs $\bar{P}_{n,-i}(\theta)=\sum_{j\neq i}P_{n,j}(\theta)/(n-1)$ for $i=1,\ldots, n$. This suggests that the an individual IRT can be locally decreasing, as long as this decrease is compensated for by increases in other IRFs at the same $\theta$ values.  
%\end{remark}

%This is similar to Assumption 5' in \cite{douglas2001asymptotic}.  Condition \ref{cond:derivbound} implies that IRFs  $P_{n,i}(\theta)$ are non-decreasing with respect to  the latent trait $\theta$. 
%\begin{condition}\label{cond:unifbd}
%There exist universal constants $ f<g \in [0,1]$ such that ${P}_{n, i}(\theta) \in[f,g]$ when $n$ is sufficiently large. Moreover, for any $\epsilon>0$, there exist constants $l_{\epsilon}$ and $u_{\epsilon}\in (0,1)$ and  $N_{\epsilon}>0$ 	such that when $n\geqslant N_{\epsilon}$, 
%\begin{align*}
%	\max_{1\leqslant i\leqslant n}\sup_{\theta\in [0,l_{\epsilon}]} (P_{n,i}(\theta)-f)\leqslant \epsilon,\hspace{1.5em}\text{and} \hspace{1.5em} 	\max_{1\leqslant i\leqslant n}\sup_{\theta\in [u_{\epsilon},1]} (g-P_{n,i}(\theta)) \leqslant \epsilon. 
%\end{align*}
%\end{condition} 

\begin{condition}\label{cond:unifbd}
%Let  $f<g \in [0,1]$  be two universal constants.
For each $i\in \{1,\ldots,n\}$, 
there exist  constants $\kappa_{n,i}<\gamma_{n,i}\in [0,1]$ such that ${P}_{n, i}(\theta) \in[\kappa_{n,i},\gamma_{n,i}]$. 
Moreover, 
for any $\epsilon > 0$, there exist constants $l_{\epsilon}$ and $u_{\epsilon}\in (0,1)$  	such that for all $(n,i)$, 
%when $n\geqslant N_{\epsilon}$, for all $i=1,\ldots, n$, 
%\begin{align*}
% 	&0\leqslant  P_{n,i}(\theta)-f  \leqslant \epsilon,\hspace{1.8em}\forall \theta\in[0,\, l_{\epsilon}],\\
% 	&0\leqslant  g-P_{n,i}(\theta)   \leqslant \epsilon,  \hspace{1.8em}\forall \theta\in[u_{\epsilon},1]. 
%\end{align*}
\begin{align*}
	\ \sup_{\theta\in [0,l_{\epsilon}]} [P_{n,i}(\theta)-\kappa_{n,i}]\leqslant \epsilon,\hspace{1.5em}\text{and} \hspace{1.5em} 	 \sup_{\theta\in [u_{\epsilon},1]} [\gamma_{n,i}-P_{n,i}(\theta)]\leqslant \epsilon. 
\end{align*}
\end{condition} 

Condition \ref{cond:unifbd} implies that for each $(n,i)$, $\lim_{\theta\downarrow 0} P_{n,i}(\theta)=\kappa_{n,i}$ and $\lim_{\theta\uparrow 1} P_{n,i}(\theta)=\gamma_{n,i}$,
where $\theta\downarrow 0$ and $\theta\uparrow 0$ represent the one-sided limits  ``from above'' and   ``from below'', respectively.  
%Assumption 5 in \cite{douglas2001asymptotic} holds with $\kappa_{n,i}=0$ and $\gamma_{n,i}=1$.
When choosing $\kappa_{n,i}=0$ and $\gamma_{n,i}=1$, Condition \ref{cond:unifbd} implies that Assumption 5 in \cite{douglas2001asymptotic} holds, i.e., the IRFs converge to 0 and 1 on the two end points of $(0,1)$, respectively. 
In contrast, Condition \ref{cond:unifbd} allows more flexible limiting values. This would  enlarge the model class to include models with guessing and missing parameters; please see more detailed discussions in Example 2 below. 

%enable the model to incorporate   guessing and missing parameters  that will be illustrated below.    
%\blue{(change up to constant so that three-parameter model can be included)}
%This is also uniform convergence. 

In summary, 
Conditions \ref{cond:unilocal}--\ref{cond:uniformdist} are the same as Assumptions 1--3 in \cite{douglas2001asymptotic}, whereas    Conditions \ref{cond:derivbound} and \ref{cond:unifbd}   considerably relax  Assumptions 4 and 5 in \cite{douglas2001asymptotic}, respectively. 
To demonstrate this, we next provide examples  with rigorous theoretical justifications.

\vspace{0.1em}
\subsection{Examples}  \label{sec:example}
Before presenting specific examples, we point out that the model of manifest distributions consists of both the IRFs and the distribution of the latent trait. In practice, a non-uniform distribution of the latent trait, e.g., the standard  normal distribution, may be more commonly used. 
As discussed as after Condition \ref{cond:uniformdist}, 
$\Theta\sim U(0,1)$ represents just one specific parameterization of the latent trait. 
Given a continuous distribution of the latent trait that is not $U(0,1)$, we can reparametrize the latent trait and the IRFs to obtain an equivalent model.  

For instance, consider an item response model with an IRF denote as $Q(\cdot)$, and it is  coupled with the latent trait $\Lambda$ following a distribution with the cumulative  distribution function $F$.
Suppose $\Lambda$ has a density function $f(\lambda)$, and  $F$ has an inverse function denoted as $F^{-1}$. 
By Proposition 3.1 in \cite{embrechts2013note}, 
 $\Theta=F(\Lambda)\sim U(0,1)$,  and $F^{-1}(\Theta)$ is a random variable with the cumulative distribution $F$. 
Then we can construct an equivalent item response model with the IRF 
\begin{align} \label{eq:transform}
	P(\theta)=Q[ F^{-1}(\theta) ]. 
\end{align}
%When the first-order derivative of $F$ exists, 
By the chain rule and the inverse function theorem in calculus,  
\begin{align*}
	P'(\theta)=\frac{Q'[ F^{-1}(\theta) ]}{f[F^{-1}(\theta)]}.
\end{align*}
%where $f(\lambda)=F'(\lambda)$ represents the density function of $\Lambda$. 
In the classical IRT models, 
it is common to assume that the latent trait 
 $\Lambda\sim N(0,1)$, i.e., the standard normal distribution.
 Then we can plug in $F(\cdot)=\Phi(\cdot)$, where $\Phi(\cdot)$ represents the cumulative distribution function of $N(0,1)$.   

%In the classical IRT models, 
%to model the latent trait following the standard normal distribution as in the classical IRT models, 
%other parameterization of 
%the latent trait  may be more commonly used. 

\paragraph{Example 1. Two-Parameter Normal Ogive Model.} Consider the normal ogive model 
where the latent trait $\Lambda\sim N(0,1)$, and  each  IRF  $Q_{n,i} (\lambda)=\Phi(a_{n,i}(\lambda-b_{n,i}))$ is   determined by two parameters $(a_{n,i},b_{n,i})$.  
% This is also the example that \cite{douglas2001asymptotic} discussed as violating Assumption 4 in 
By \eqref{eq:transform},  we have an equivalent model with $\Theta\sim U(0,1)$ and the IRF 
\begin{align}\label{eq:normalogivemodel}
	P_{n,i}(\theta) =&~ \Phi(a_{n,i}[\Phi^{-1}(\theta)-b_{n,i} ]), \hspace{1.5em}\text{and}\hspace{1.5em}
%\end{align*}
%and thus, 
%\begin{align*}
	P^{\prime}_{n,i}(\theta)=\frac{a_{n,i} \phi(a_{n,i}[\Phi^{-1}(\theta)-b_{n,i}]) }{ \phi[\Phi^{-1}(\theta)]}, 
\end{align}
where $\phi(x)$ denotes the  density function of $N(0,1)$. 
 \cite{douglas2001asymptotic} has pointed out that this model would violate  their Assumption 4. 
In Proposition \ref{prop:normalogive} below,
we formally demonstrate that $P'_{n,i}(\theta)$ cannot be uniformly bounded over the entire interval $(0,1)$,  
Furthermore, we prove that Conditions \ref{cond:derivbound} and \ref{cond:unifbd} are  satisfied as long as the parameters $(a_{n,i}, b_{n,i})$ are uniformly bounded.

%We can formally prove that in this case, $P'_{n,i}$ cannot satisfy 
%%\begin{align*}
%%	\lim_{\theta \to 0} P'_{n,i}(\theta)\to 0
%%\end{align*}
%%so that 
%the uniform boundedness of $P'(\theta)$ in Assumption 4 of \cite{douglas2001asymptotic}. % cannot be satisfied. 
%On the other hand, we can prove that when the parameters $(a_{n,i}, b_{n,i})$ are uniformly bounded, 

\begin{proposition}\label{prop:normalogive} Suppose   IRFs in the sequence $\{Q_{n,i}(\lambda)\}$  follow the normal ogive models  with parameters $(a_{n,i}, b_{n,i})$ and the latent trait $\Lambda\sim N(0,1)$. Let  $\{P_{n,i}(\theta)\}$ denote the corresponding transformed IRFs  following \eqref{eq:normalogivemodel}  with the latent trait $\Theta=\Phi(\Lambda)\sim U(0,1)$.  
\begin{itemize}
	\item[(i)] When $|a_{n,i}|\neq 1$ or $b_{n,i}\neq 0$,  $\lim_{\theta \downarrow 0 }P'_{n,i}(\theta)$ and $\lim_{\theta \uparrow 1 }P'_{n,i}(\theta)$ are either 0 or $+\infty$. 
\item[(ii)] Assume there exist    constants $C_a, C_b>0$ independent with $(n,i)$ such that $a_{n,i}\in  [1/C_a, C_a]$ and $|b_{n,i}|\in [1/C_b, C_b]$ for all $(n,i)$. Then the transformed IRFs $\{P_{n,i}(\theta)\}$ satisfy  Conditions \ref{cond:derivbound} and \ref{cond:unifbd} with $\kappa_{n,i}=0$ and $\gamma_{n,i}=1$.   
\end{itemize}
\end{proposition}

\begin{figure}[!htbp]
\centering
\begin{subfigure}{0.32\textwidth}
\caption{$P_{n,i}(\theta)$ for $\theta \in (0, 1)$}
\includegraphics[width=\textwidth]{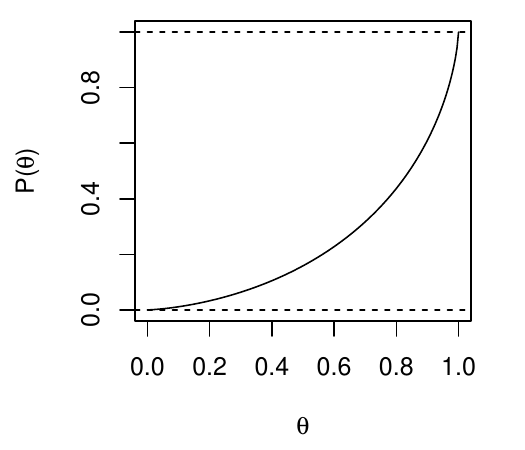}
\end{subfigure}	
\begin{subfigure}{0.32\textwidth}
\caption{$P_{n,i}'(\theta)$ for  $\theta \downarrow 0$}
\includegraphics[width=\textwidth]{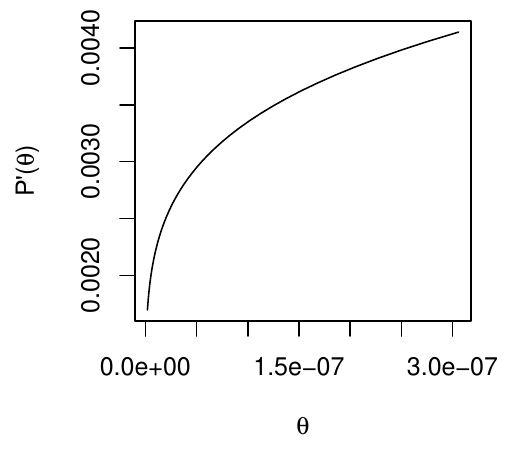}
\end{subfigure}	
\begin{subfigure}{0.32\textwidth}
\caption{$P_{n,i}'(\theta)$ for  $\theta \uparrow 1$} 
\includegraphics[width=\textwidth]{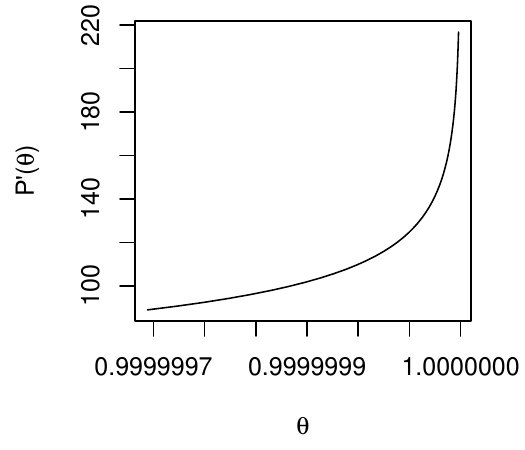}
\end{subfigure}	
\caption{$P_{n,i}(\theta)$ and $P'_{n,i}(\theta)$ for the normal ogive model in \eqref{eq:normalogivemodel} with $a_{n,i}=b_{n,i}=1$.}\label{fig:normalogive}
\end{figure}

We visually illustrate Proposition \ref{prop:normalogive} by plotting  $P_{n,i}(\theta)$ over $(0,1)$ and and $P'_{n,i}(\theta)$ when $\theta$ is close to 0 and 1, respectively. 
Figure \ref{fig:normalogive} suggests that $P'_{n,i}(\theta)$ approaches infinitesimal proximity to 0 and $\infty$ as $\theta$ converges to 0 and 1, respectively, which is consistent with Proposition \ref{prop:normalogive} (a).  
Nevertheless, when $\theta$ is bounded away from 0 and 1, $P'(\theta)$ is finite and strictly positive.  
The above discussions focus on  two-parameter normal ogive model, chosen for simplicity and to align with the discussions in \cite{douglas2001asymptotic}. 
Similar conclusions can also be established beyond this model class, even with additional  parameters in the model. As an example, we  next examine the four-parameter logistic model.

\paragraph{Example 2. Four-Parameter Logistic Model.}
Typically, the 
four-parameter logistic  (4PL) model  assumes that the latent trait $\Lambda \sim N(0,1)$ and the IRF $Q_{n,i}(\lambda)=c_{n,i}+(d_{n,i}-c_{n,i})g[a_{n,i}(\lambda-b_{n,i})]$ depends on four parameters $(a_{n,i},b_{n,i},c_{n,i},d_{n,i})$ and $g(x)=e^x/(1+e^x)$. 
When we consider the reparametrized latent trait $\Theta=\Phi^{-1}(\Lambda)\sim U(0,1)$, by \eqref{eq:transform}, the equivalently transformed IRF is 
\begin{align}\label{eq:transicc4pl}
P_{n,i}\left(\theta\right)= c_{n,i}+\left(d_{n,i}-c_{n,i}\right) g\left[a_{n,i}\left(\Phi^{-1}(\theta)-b_{n,i}\right)\right].  
\end{align}
%In this case, $f=c_i$ and $g=d_i$.  
This formulation can cover the Rasch model, 2PL, and 3PL  models as   special cases by setting some parameters to 0. 
In \eqref{eq:transicc4pl}, $\lim_{\theta\downarrow 0} P_{n,i}(\theta)=c_{n,i}$,	$\lim_{\theta\uparrow 1} P_{n,i}(\theta)=d_{n,i}$, and 
\begin{align*}
	P'_{n,i}(\theta)=(d_{n,i}-c_{n,i})a_{n,i}\frac{ g'\left[a_{n,i}\left(\Phi^{-1}(\theta)-b_{n,i}\right)\right]}{\phi[\Phi^{-1}(\theta)]}. 
\end{align*}
In Proposition \ref{prop:4pllogisticmodel} below, 
we formally show that $P'_{n,i}(\theta)$ is unbounded over the entire interval $(0,1)$ and prove that Conditions \ref{cond:derivbound} and \ref{cond:unifbd} can be  satisfied.  
\begin{proposition}\label{prop:4pllogisticmodel}
Suppose   IRFs in the sequence $\{Q_{n,i}(\lambda)\}$  follow the 4PL models  with the latent trait $\Lambda\sim N(0,1)$ and the  parameters $(a_{n,i}, b_{n,i}, c_{n,i}, d_{n,i})$ satisfying   $c_{n,i}<d_{n,i}\in [0,1]$. 
 Let  $\{P_{n,i}(\theta)\}$ denote the corresponding transformed IRFs  following \eqref{eq:transicc4pl}  with the latent trait $\Theta=\Phi(\Lambda)\sim U(0,1)$.  
\begin{itemize}
	\item[(i)] When $a_{n,i} \neq 0$,  $\lim_{\theta \downarrow 0 }P'_{n,i}(\theta)=\lim_{\theta \uparrow 1 }P'_{n,i}(\theta) = +\infty$. 
\item[(ii)] Assume there exist    constants $C_a, C_b, C_{c,d}>0$ independent with $(n,i)$  such that $a_{n,i}\in  [1/C_a, C_a]$, $|b_{n,i}| \in [1/C_b, C_b]$, and $d_{n,i}-c_{n,i}\in [1/C_{c,d}, C_{c,d}]$ for all $(n,i)$. Then the transformed IRFs $\{P_{n,i}(\theta)\}$ satisfy  Conditions \ref{cond:derivbound} and \ref{cond:unifbd} with $\kappa_{n,i}=c_{n,i}$ and $\gamma_{n,i}=d_{n,i}$.  
\end{itemize}
	
\end{proposition}

\begin{figure}[!htbp]
\centering 
\begin{subfigure}{0.32\textwidth}
\caption{$P_{n,i}(\theta)$ for $\theta \in (0, 1)$}
\includegraphics[width=\textwidth]{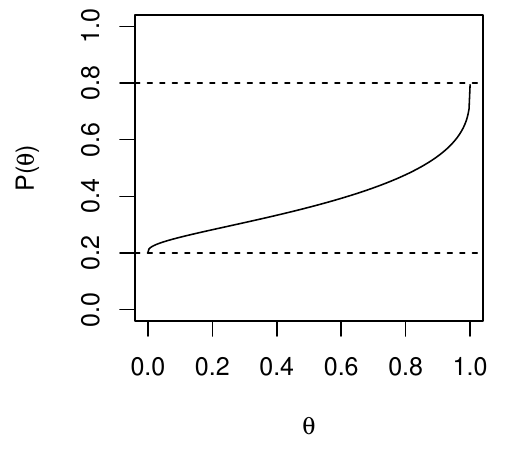}
\end{subfigure}	
\begin{subfigure}{0.32\textwidth}
\caption{$P_{n,i}'(\theta)$ for  $\theta \downarrow 0$}
\includegraphics[width=\textwidth]{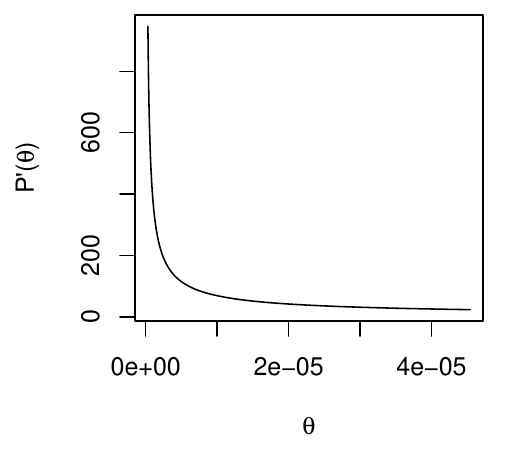}
\end{subfigure}	
\begin{subfigure}{0.32\textwidth}
\caption{$P_{n,i}'(\theta)$ for  $\theta \uparrow 1$} 
\includegraphics[width=\textwidth]{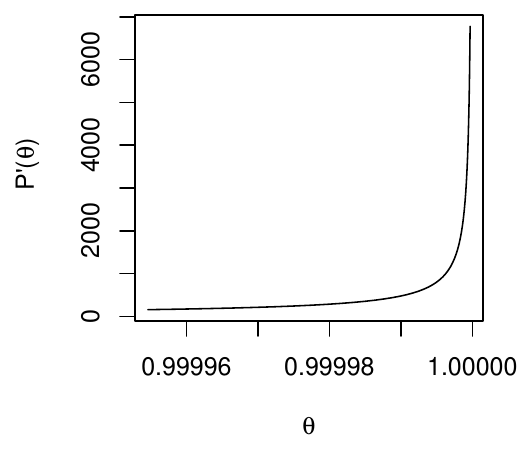}
\end{subfigure}	
\caption{$P_{n,i}(\theta)$ and $P'_{n,i}(\theta)$ for the 4PL  in \eqref{eq:transicc4pl} with $a_{n,i}=b_{n,i}=1$, $c_{n,i}=0.2$, and $d_{n,i}=0.8$.}\label{fig:4plmodel} 
\end{figure}

We visually illustrate Proposition \ref{prop:4pllogisticmodel} by plotting  $P_{n,i}(\theta)$ over $(0,1)$ and $P'_{n,i}(\theta)$ when $\theta$ is close to 0 and 1, respectively. 
Figure \ref{fig:4plmodel} (a) shows that  $P_{n,i}(\theta)$ converges to 0.2 and 0.8 at the two ends points of (0,1), respectively. Therefore    Assumption 5 in \cite{douglas2001asymptotic} is violated in this case. 
Moreover, Figure \ref{fig:4plmodel} (b)--(c) suggests that $P'_{n,i}(\theta)$ diverges to $\infty$ as $\theta$ converges to 0 and 1, which aligns with Proposition \ref{prop:4pllogisticmodel} (a) and shows that Assumption 4 in \cite{douglas1997joint} is violated.  

%However, when $\theta$ is bounded away from 0 and 1, $P'(\theta)$ is finite and strictly positive. 

%\vspace{1.5em}

\section{Results}\label{sec:results}

We next present Theorem \ref{thm:asymident} establishing the asymptotic equivalence between   two sequences of IRFs under relaxed assumptions. %, $\{\mathcal{P}_k, \mathcal{P}_{k+1},\ldots\}$ and $\{\mathcal{P}_k^*, \mathcal{P}_{k+1}^*,\ldots\}$. 

\begin{theorem}\label{thm:asymident}
For any two sequences of IRFs $\{\mathcal{P}_k, \mathcal{P}_{k+1},\ldots\}$ and $\{\mathcal{P}_k^*, \mathcal{P}_{k+1}^*,\ldots\}$,
given the sequence of the manifest distribution $\{\mathcal{F}_k, \mathcal{F}_{k+1},\ldots,\}$, under Conditions \ref{cond:unilocal}--\ref{cond:unifbd}, 
\begin{align*}
\lim_{n\to \infty}	\max_{1 \leq i \leq n} \sup_{\theta \in(0, 1)}\left|P_{n, i}(\theta)-P_{n, i}^*(\theta)\right| = 0.
\end{align*}	 
\end{theorem}

It is  important to note that the asymptotic equivalence in  Theorem \ref{thm:asymident} holds over the entire interval $(0,1)$ even though that Condition \ref{cond:derivbound} only assumes bounded derivatives on compact subsets of $(0,1)$. 
Intuitively, this is achievable because the value of  $P_{n,i}(\theta)$ has limited variation near the two end points of $(0,1)$, as specified by Condition \ref{cond:unifbd}. 
However, we emphasize that the proof is not a simple application of Condition \ref{cond:unifbd}. 
The relaxations introduced by Conditions \ref{cond:derivbound}--\ref{cond:unifbd}  necessitate the development of novel theoretical techniques. 
Due to the relaxations introduced by Conditions \ref{cond:derivbound}--\ref{cond:unifbd}, many arguments in \cite{douglas2001asymptotic} are no longer applicable.  
We need to develop novel theoretical techniques  that adapt to different ranges of $\theta$.  
Please also see the detailed proofs and more technical discussions in Remark \ref{rm:novelty}. 

%Please also see the detailed proofs below. 

%\begin{remark}
%	  \cite{douglas1997joint} studied the joint consistency under an assumption similar to Condition \ref{cond:derivbound}. Nevertheless, we note that that the results in \cite{douglas1997joint} is derived on a fixed compact interval subset of $(0,1)$. Using the techniques in this paper, we may also extend the results in \cite{douglas1997joint} to the whole interval $(0,1)$. 
%\end{remark}

\vspace{1.5em}
\textit{Proof.} To prove Theorem \ref{thm:asymident},
it suffices to prove that for any $\epsilon>0$, 
there exists $N_{\epsilon}>0$ such that when $n\geqslant N_{\epsilon}$, 
$
	\max_{1 \leq i \leq n} \sup_{\theta \in(0, 1)}|P_{n, i}(\theta)-P_{n, i}^*(\theta)| \leqslant C\epsilon
$ 
where $C>0$ is a universal constant. 
Given any $\epsilon>0$, 
let $l_{\epsilon}$ and $u_{\epsilon}$ be defined as in Condition \ref{cond:unifbd}. Then we have 
$\max_{1\leqslant i\leqslant n}\sup_{\theta \in(0, 1)}|P_{n, i}(\theta)-P_{n, i}^*(\theta)|\leqslant \max\{B_1,B_2,B_3\}$, where
\begin{align*}
&B_1=\max_{1\leqslant i\leqslant n}	\sup_{\theta \in(0, l_{\epsilon}]}\left|P_{n, i}(\theta)-P_{n, i}^*(\theta)\right|,\hspace{2em} B_2=	\max_{1\leqslant i\leqslant n}\sup_{\theta \in [u_{\epsilon},1)}\left|P_{n, i}(\theta)-P_{n, i}^*(\theta)\right|\\
&B_3=	\max_{1\leqslant i\leqslant n}\sup_{\theta \in(l_{\epsilon},u_{\epsilon})}\left|P_{n, i}(\theta)-P_{n, i}^*(\theta)\right|. 
\end{align*}
By Condition \ref{cond:unifbd}, 
\begin{align*}
	B_1 =&~ \max_{1\leqslant i\leqslant n}\sup_{\theta \in(0, l_{\epsilon}]}\left|(P_{n, i}(\theta)-\kappa_{n,i})-(P_{n, i}^*(\theta)-\kappa_{n,i})\right|\notag\\
	\leqslant &~\max_{1\leqslant i\leqslant n}\sup_{\theta \in(0, l_{\epsilon}]}\left\{\left| P_{n, i}(\theta)-\gamma_{n,i}  \right|+\left| P_{n, i}^*(\theta)-\gamma_{n,i} \right|\right\}
	\leqslant  2\epsilon,\\ 
%\end{align*}
%Similarly, by Condition \ref{cond:unifbd},
%\begin{align*}
	B_2\leqslant &~\max_{1\leqslant i\leqslant n}\sup_{\theta \in [ u_{\epsilon},1)}\left\{\left| P_{n, i}(\theta)-\gamma_{n,i}   \right|+\left| P_{n, i}^*(\theta)-\gamma_{n,i} \right|\right\}
	\leqslant  2\epsilon. 
\end{align*}
%\begin{align*}
%	&~\sup_{\theta \in(0, 1)}\left|P_{n, i}(\theta)-P_{n, i}^*(\theta)\right| \notag\\
%\leqslant &~\sup_{\theta \in(0, l_{\epsilon}]}\left|P_{n, i}(\theta)-P_{n, i}^*(\theta)\right| \notag\\
%\end{align*}
%For any $\epsilon>0$,
%by Condition \ref{cond:unifbd}, 
%there exists a constant $l_{\epsilon}\in (0,1)$ that is independent with $(n,i)$, such that
%% \blue{there exists $l_{\epsilon}\in (0,1)$ (what assumption?)} such that 
%\begin{align*}
%	\sup_{\theta\in (0,l_{\epsilon})}P_{n,i}(\theta)\leqslant \epsilon,\hspace{1.5em}\text{and} \hspace{1.5em}\sup_{\theta\in (u_{\epsilon},1)}\{1-P_{n,i}(\theta)\}\leqslant \epsilon. \end{align*}
%\begin{align*}
%	\sup_{\theta\in (0, l_{\epsilon})} P_{n,i}(\theta)\leqslant P_{n,i}(l_{\epsilon})\leqslant \epsilon. 
%\end{align*}	
%Thus
%\begin{align*}
%	\max_{1\leqslant i\leqslant n} 	\sup_{\theta\in (0, l_{\epsilon})} |P_{n, i}(\theta)-P_{n, i}^*(\theta)|\leqslant 2\epsilon. 
%\end{align*}
%Similarly,  there exists $u_{\epsilon}\in (0,1)$ such that 
%\begin{align*}
%	\sup_{\theta\in (u_{\epsilon},1)} 1-P_{n,i}(\theta)\leqslant 1-P_{n,i}(u_{\epsilon})\leqslant \epsilon. 
%\end{align*}	
 We next establish Theorem \ref{thm:asymepsilonbd}    showing that there exists $N_{\epsilon}>0$ such that when $n\geqslant N_{\epsilon}$, $B_3\leqslant \epsilon$. Then the proof of Theorem \ref{thm:asymident} is finished.  
%Moreover, by Proposition  \ref{thm:asymepsilonbd}, there exists $N_{\epsilon}>0$ such that when $n\geqslant N_{\epsilon}$, $B_3\leqslant \epsilon$. 

%\begin{align*}
%	\max_{1\leqslant i\leqslant n} 	\sup_{\theta\in (l_{\epsilon}, u_{\epsilon})} |P_{n, i}(\theta)-P_{n, i}^*(\theta)|\leqslant B_{l_{\epsilon}, u_{\epsilon}} n^{-\alpha}\leqslant \epsilon.  
%\end{align*}
%There exists $n_{0,\epsilon}$ such that when  $n\geqslant  n_{0,\epsilon}$, we have 
%\begin{align*}
%	\max_{1\leqslant i\leqslant n} 	\sup_{\theta\in (l_{\epsilon}, u_{\epsilon})} |P_{n, i}(\theta)-P_{n, i}^*(\theta)|\leqslant 2\epsilon. 
%\end{align*}

\vspace{2em}

%We next present Theorem \ref{thm:asymepsilonbd}, whose proof relies on three lemmas stated and proved in the Appendix section. 

\vspace{0.5em}
\begin{theorem}\label{thm:asymepsilonbd}
Assume Conditions \ref{cond:unilocal}--\ref{cond:unifbd}. 
For any given $\epsilon>0$ and $\alpha <\beta \in (0,1)$, there exists $N_{\epsilon,\alpha,\beta}$ such that when $n\geqslant N_{\epsilon,\alpha,\beta}$,
\begin{align*}
	\max_{1 \leq i \leq n} \sup_{\theta \in(\alpha,\beta)}\left|P_{n, i}(\theta)-P_{n, i}^*(\theta)\right|<\epsilon. 
\end{align*}
\end{theorem} 
%\begin{align*}
%\end{align*}	
%where $C_{\alpha\beta}(\alpha)$ represents a constant that only depends on $\alpha,\beta,$ and $\alpha$. \blue{(I might change this to: For any $\epsilon>0$, there exists $N_{\epsilon}$ such that when $n\geqslant N_{\epsilon}$, )}
%we need to develop sophisticated investigations considering different $\theta$ ranges. s 

\vspace{1em}

\textit{Proof.} 
Consider an arbitrary item $i\in \{1,\ldots, n\}$ and $\theta\in (\alpha,\beta)$. 
For each integer $k$ such that $\alpha <k/(n-1)< \beta$, define $\theta_k$ and $\theta_k^*$ to satisfy 
\begin{align}\label{eq:defthetak}
	\bar{P}_{n,-i}\left(\theta_k\right)=\bar{P}_{n,-i}^*\left(\theta_k^*\right)=k /(n-1),
\end{align}
where we define the functions  
\begin{align*}
	\bar{P}_{n,-i}(\theta)=\sum_{j\neq i} P_{n,i}(\theta)/(n-1),\quad\quad\text{and}\quad\quad \bar{P}_{n,-i}^*(\theta)=\sum_{j\neq i} P_{n,i}^*(\theta)/(n-1), 
\end{align*}
 representing the means of the IRFs of their respective sequences. 
 Note that $\theta_k$ and $\theta_k^*$  depend on $n$ and $i$, but this is suppressed in the notation for simplicity of presentation.  
For any given $\theta\in (\alpha,\beta)$,  
select the integer $k$ such that $\theta_k\in (\alpha,\beta)$ and $|\theta-\theta_k|$ is minimized. Then
\begin{align*}
	 &~|P_{n,i}(\theta)-P_{n,i}^*(\theta)| \notag\\
\leqslant &~|P_{n,i}(\theta)-P_{n,i}(\theta_k)|+|P_{n,i}(\theta_k)-P_{n,i}^*(\theta_k^*)|  +|P_{n,i}^*(\theta_k^*)-P_{n,i}^*(\theta)|.  \notag
\end{align*}
The following proof consists of two main steps showing that
\begin{align*}
	\text{Step 1: }&\quad \max\{|P_{n,i}(\theta)-P_{n,i}(\theta_k)|,\ |P_{n,i}^*(\theta_k^*)-P_{n,i}^*(\theta)|\}\leqslant  \frac{4M_{\alpha\beta}}{m_{\alpha\beta}n}, \\
	\text{Step 2: }&\quad 	|P_{n,i}(\theta_k)-P_{n,i}^*(\theta_k^*)|  \leqslant  \epsilon,
\end{align*}
respectively, where $m_{\alpha\beta}$ and $M_{\alpha\beta}$ are constants specified as in Condition \ref{cond:derivbound}. 
It is worth mentioning that $M_{\alpha\beta}$ and $m_{\alpha\beta}$ are constants that depend on $(\alpha,\beta)$ but are independent with $(n,i)$.

% by Condition  \ref{cond:derivbound}. 

\paragraph{Step 1.} 
%Let $ F_{n,-i}$ denote the cumulative distribution function of $\bar{Y}_{n,-i}$. Then
%\begin{align*}
%	 F_{n,-i}\left(\bar{P}_{n,-i}\left(\theta_k\right)\right)=F_{n,-i}\left(\bar{P}_{n,-i}^*\left(\theta_k^*\right)\right)=F_{n,-i}(k /(n-1))
%\end{align*}
%Derive upper bound of $|\theta-\theta_k|$ and then we can bound $|P_{n,i}(\theta)-P_{n,i}(\theta_k)|$ and $|P_{n,i}^*(\theta_k^*)-P_{n,i}^*(\theta)|$. 
%The upper bound of $|\theta-\theta_k|$ would require derivative lower bound. (but can we relax it to fixed region?) 

As $\theta$ and $\theta_k \in (\alpha,\beta)$,  
\begin{align} \label{eq:step1bd1}
	|P_{n,i}(\theta)-P_{n,i}(\theta_k)|\leqslant &~\sup_{\eta\in (\alpha,\beta)} |P'_{n,i}(\eta)|\times |\theta-\theta_k| \leqslant M_{\alpha\beta} \times |\theta-\theta_k|. 
\end{align} 
%where $M_{\alpha\beta}$ is a constant that depends on $(\alpha,\beta)$ but is independent with $(n,i)$ by Condition  \ref{cond:derivbound}. 
%Moreover, by Condition  \ref{cond:bound}, 
%we know $\inf_{\theta\in (\alpha,\beta)} |P'_{n,i}(\theta)|>m_{\alpha\beta}$. 
%\begin{align*}
%	|\bar{P}_{n,-i}(\theta_k)-\bar{P}_{n,-i}(\theta)|
%\end{align*}
%Given $\theta_k$, let $\theta_{k-1}$ and $\theta_{k+1}$ to be the second  closet values to $\theta$ that are smaller and larger than $\theta_k$, respectively, i.e., $\theta_{k-1}<\theta_k$ and $$ 
By the definition in \eqref{eq:defthetak} and Condition \ref{cond:unifbd},
there exist $(l_{\alpha,\beta}, u_{\alpha,\beta})$ independent with $(n,i)$ such that $0<l_{\alpha,\beta}<  \theta_{k-1}, \theta_{k+1}  < u_{\alpha,\beta}<1$ when $n$ is sufficiently large.   
%$\theta_{k-1}, \theta_{k+1} $ 
%By Condition \ref{cond:unifbd}, 
By Condition \ref{cond:derivbound}, 
$P_{n,i}(\theta)$ and  $\bar{P}_{n,-i}(\theta)$  are non-decreasing functions with respect to $\theta$ on a fixed interval $(l_{\alpha,\beta}, u_{\alpha,\beta})$.   
%By the \blue{non-decreasing assumption} of $P_{n,i}$ (and thus $\bar{P}_{n,-i}$ also),
Therefore,  we know  $\theta_{k-1}<\theta_k < \theta_{k+1}$. 
% $P_{n,i}(\theta_k)\leqslant 
%P_{n,i}(\lambda_{k+})\leqslant P_{n,i}(\theta_{k+1})$. 
%Define $\lambda_{k-}=\max\{\theta_{k-1}, a \}$ and $\lambda_{k+}=\min\{\theta_{k+1}, b\}$. 
%As $k$ is the integer that minimizes $|\theta-\theta_k|$, and $\theta,\theta_k\in (\alpha,\beta)$, 
% $\theta_k$ is the closet one to $\theta$, and $\theta,\theta_k\in (\alpha,\beta)$, 
%By $\theta,\theta_k\in (\alpha,\beta)$, 
Moreover, by $\theta,\theta_k\in (\alpha,\beta)$, we have $$\alpha \leqslant \lambda_{k-}<\theta_k < \lambda_{k+} \leqslant \beta,$$  where we define
 $\lambda_{k-}=\max\{\theta_{k-1}, \alpha \}$ and $\lambda_{k+}=\min\{\theta_{k+1},  \beta\}$. 
% we have $a\leqslant \lambda_{k-}<\theta_k < \lambda_{k+} \leqslant b$. Therefore, 
As $k$ is the integer that minimizes $|\theta-\theta_k|$, 
\begin{align}
	|\theta-\theta_k|\leqslant &~ |\lambda_{k-}-\lambda_{k+}|  
\leqslant  \frac{|\bar{P}_{n,-i}(\lambda_{k-})-\bar{P}_{n,-i}(\lambda_{k+})|}{\inf_{\eta\in [\alpha,\beta]}|\bar{P}_{n,-i}'(\eta) |}. \label{eq:thetathetakbd}
\end{align}
By Condition  \ref{cond:derivbound}, there exists a constant $m_{\alpha\beta}>0$ independent with $(n,i)$ such that   $\inf_{\eta\in (\alpha,\beta)} |P'_{n,i}(\eta)|>m_{\alpha\beta}$. Thus,  
\begin{align*}
	\eqref{eq:thetathetakbd} \leqslant  &~\frac{1}{m_{\alpha\beta}}|\bar{P}_{n,-i}(\lambda_{k-})-\bar{P}_{n,-i}(\lambda_{k+})| \notag\\
	\leqslant &~\frac{1}{m_{\alpha\beta}}|\bar{P}_{n,-i}(\theta_{k-1})-\bar{P}_{n,-i}(\theta_{k+1})| = \frac{2}{m_{\alpha\beta}(n-1)}\leqslant \frac{4}{m_{\alpha\beta}n}
\end{align*}
where the second inequality is obtained by the monotonicity of $\bar{P}_{n,-i}(\theta)$ under Condition \ref{cond:unifbd}. 
%\begin{align*}
%|\theta-\theta_k|\leqslant &~ |\lambda_{k-}-\lambda_{k+}|  
%\leqslant  \frac{|\bar{P}_{n,-i}(\lambda_{k-})-\bar{P}_{n,-i}(\lambda_{k+})|}{\inf_{\eta\in [\alpha,\beta]}|\bar{P}_{n,-i}'(\eta) |} \notag\\
%\leqslant &~ \frac{|\bar{P}_{n,-i}(\theta_{k-1})-\bar{P}_{n,-i}(\theta_{k+1})|}{\inf_{\eta\in [\alpha,\beta]}|\bar{P}_{n,-i}'(\eta) |}\leqslant \frac{2}{m_{\alpha\beta}(n-1)}, 
%\end{align*}  
%\begin{align*}
%|\theta-\theta_k|\leqslant &~ \max\{|\lambda_{k-}-\theta_k|, |\lambda_{k+}-\theta_k| \} \notag\\
%\leqslant &~ \frac{\max\{|\bar{P}_{n,-i}(\lambda_{k-})-\bar{P}_{n,-i}(\theta_k)|,\ |\bar{P}_{n,-i}(\lambda_{k+})-\bar{P}_{n,-i}(\theta_k)| \}}{\inf_{\eta\in [\alpha,\beta]}\bar{P}_{n,-i}'(\eta) } \notag\\
%\leqslant &~   \frac{1}{m_{\alpha\beta}(n-1)}. 
%\end{align*}    
%\begin{align*}
%	|\theta-\theta_k|\leqslant \max_{j\in \{k-1,k+1\}}|\theta_j-\theta_k|\leqslant \max_{j\in \{k-1,k+1\}}\frac{|\bar{P}_{n,-i}(\theta_j)-\bar{P}_{n,-i}(\theta_k)|}{\inf_{\eta\in (\alpha,\beta)} \bar{P}'_{n,-i}(\eta) }\leqslant \frac{1}{m_{\alpha\beta}(n-1)}. 
%\end{align*} 
%where in the last equality,  we use that by Condition  \ref{cond:derivbound}, there exists a constant $m_{\alpha\beta}>0$ independent with $(n,i)$ such that   $\inf_{\eta\in (\alpha,\beta)} |P'_{n,i}(\eta)|>m_{\alpha\beta}$. 
In summary, we have $\eqref{eq:step1bd1} \leqslant {4M_{\alpha\beta}}/{(m_{\alpha\beta}n)}.$
The same upper bound can be obtained for $|P_{n,i}^*(\theta_k^*)-P_{n,i}^*(\theta)|$ following a similar analysis, and thus Step 1 is proved.

\paragraph{Step 2.} 
Define $\bar{Y}_{n,-i} = \sum_{j\neq i}Y_{n,i}/(n-1) $ and the event $\mathcal{E}_{n,k}=\{\bar{Y}_{n,-i}=k /(n-1)\}$. 
Let $\delta\in (0,1/2)$ be a fixed small number. Define an interval $I_{\delta}=(\delta, 1-\delta)$. 
%We derive the upper bound using  
Then we have 
\begin{align*}
&~|P_{n,i}(\theta_k)-P_{n,i}^*(\theta_k^*)|\notag\\
\leqslant &~ |P_{n,i}(\theta_k)-P(Y_{n,i}=1, \, \Theta\in I_{\delta}\mid \mathcal{E}_{n,k})
|+|P(Y_{n,i}=1, \, \Theta\in I_{\delta}\mid \mathcal{E}_{n,k})
-P_{n,i}^*(\theta_k^*)|\notag \\
\leqslant &~ A_1+A_2+A_1^*+A_2^*, 
\end{align*}  
%We have
%\begin{align*}
%&~ \left| P_{n, i}\left(\theta_k\right)-P\left(Y_{n, i}=1,\, \Theta\in I_{\delta} \mid \mathcal{E}_{n,k}\right) \right|\leqslant A_1+A_2	
%\end{align*}
where we define
\begin{align*}
	A_1=&~ \left| P_{n, i}\left(\theta_k\right)\big\{1-P(\Theta\in I_{\delta}\mid \mathcal{E}_{n,k})\big\} \right|   \\
	A_2=&~\left|P_{n, i}\left(\theta_k\right)P\left(\Theta\in I_{\delta}\mid \mathcal{E}_{n,k}\right)-P\left(Y_{n, i}=1, \Theta\in I_{\delta}\mid \mathcal{E}_{n,k} \right)\right|\\
	A_1^*=&~ \left| P_{n, i}^*\left(\theta_k^*\right)\big\{1-P(\Theta\in I_{\delta}\mid \mathcal{E}_{n,k})\big\} \right| \\
	A_2^*=&~\left|P_{n, i}^*\left(\theta_k^*\right)P\left(\Theta\in I_{\delta}\mid \mathcal{E}_{n,k}\right)-P\left(Y_{n, i}=1, \Theta\in I_{\delta}\mid \mathcal{E}_{n,k} \right)\right|. 
\end{align*}
We point out that   $(A_1,A_1^*,A_2,A_2^*)$   depend on $(n,i,k,\delta)$, but this is suppressed in the notation for simplicity. 
By Lemma \ref{lm:fixedrange} in the appendix, 
\begin{align}\label{eq:a1bddelta}
A_1+A_1^*\leqslant \frac{2 n\delta\exp(-n\tilde{C}_{\alpha \beta,1} )}{\tilde{C}_{\alpha\beta,2}}
%\frac{4n\delta\exp[-(n-1)\tilde{C}_{\alpha}^2 /2  ]}{\tilde{C}_{\alpha\beta}}	
\end{align}
As the exponential term converges to 0 faster than polynomial $n^{-r}$ for any $r>0$. Therefore, for any $\epsilon>0$, there exists $N_{\epsilon,\alpha,\beta,\delta}$ such that for $n\geqslant N_{\epsilon,\alpha,\beta,\delta}$, $A_1+A_1^*\leqslant \epsilon/ 2$. 
%Therefore, when $n$ is sufficiently large, we have $\eqref{eq:a1bddelta}\leqslant \epsilon/2$. 
%when $\delta$ is sufficiently small, there exists a constant $C_{\alpha,\beta,2}$ independent with $(n,i)$ such that   
%\begin{align}\label{eq:a1bddelta}
%	A_1+A_1^*\leqslant {C_{\alpha,\beta,2}}\, {\max_{j\neq i} P_{n,j}(\delta) \delta}. 
%\end{align} 
%As IRFs are increasing functions by Condition \ref{cond:derivbound}, 
%By Conditions \ref{cond:derivbound} and \ref{cond:unifbd}, 
%we can 
%choose $\delta_{\epsilon}$ that is sufficiently small and independent with $(n,i)$,  such that 
%when $\delta\leqslant \delta_{\epsilon}$, $\eqref{eq:a1bddelta}\leqslant \epsilon/2$. 

%\begin{lemma}\label{lm:diffbdinterval}
%For any fixed $\epsilon>0$ and $\delta>0$, there exists $N_{\epsilon,\alpha,\beta,\delta}$ that is independent with $i$ such that  when $n\geqslant N_{\epsilon,\alpha,\beta,\delta}$, 
%\begin{align*}
%	\left|P_{n, i}\left(\theta_k\right)-P\left(Y_{n, i}=1, \Theta\in I_{\delta} \mid \mathcal{E}_{n,k}\right)\right| \leqslant \epsilon. 
%\end{align*} 
%\end{lemma}

We next prove that for any fixed $\epsilon>0$ and $\delta>0$, when $n$ is sufficiently large, 
%there exists $N_{\epsilon,\alpha,\beta,\delta}$ that is independent with $i$ such that  when $n\geqslant N_{\epsilon,\alpha,\beta,\delta}$,  
\begin{align}\label{eq:a2a2starbd}
	A_2+A_2^*\leqslant \epsilon/2. 
\end{align}
%\begin{align*}
%	\left|P_{n, i}\left(\theta_k\right)-P\left(Y_{n, i}=1, \Theta\in I_{\delta} \mid \mathcal{E}_{n,k}\right)\right| \leqslant \epsilon. 
%\end{align*} 
%
%Moreover, by Lemma \ref{lm:diffbdinterval},
%there exists $N_{\epsilon,\alpha,\beta}>0$ that is independent with $i$ such that  
%when $n\geqslant N_{\epsilon,\alpha,\beta}$,  $A_2+A_2^*\leqslant \epsilon/2$. 
%\subsection*{Proof of Lemma \ref{lm:diffbdinterval}}
%\blue{(maybe move this to theorem)}
In particular, define an interval $I_{k,\eta}=(\theta_k-n^{-\eta}, \theta_k+n^{-\eta})$ with $0<\eta<1/2$. Then 
\begin{align*}
A_2= &~\left|P_{n, i}\left(\theta_k\right)P\left(\Theta\in I_{\delta}\mid \mathcal{E}_{n,k}\right)-P\left(Y_{n, i}=1, \Theta\in I_{\delta} \mid \mathcal{E}_{n,k}\right)\right|\notag\\
 \leqslant &~\big|P_{n, i}\left(\theta_k\right)[P(\Theta\in I_{\delta}\mid \mathcal{E}_{n,k})- P(\Theta\in I_{\delta}\cap I_{k,\eta}\mid \mathcal{E}_{n,k} )]\big| \notag\\
&~+\big|P_{n,i}(\theta_k)P\left(\Theta\in I_{\delta}\cap I_{k,\eta}\mid \mathcal{E}_{n,k} \right) - P(Y_{n,i}=1, \Theta\in I_{\delta}\cap I_{k,\eta}\mid \mathcal{E}_{n,k} )  \big| \notag\\
&~+\big| P(Y_{n,i}=1, \Theta\in I_{\delta}\cap I_{k,\eta}\mid \mathcal{E}_{n,k} ) -P(Y_{n,i}=1, \Theta\in I_{\delta} \mid \mathcal{E}_{n,k} )  \big| \notag\\
\leqslant &~ A_{21}+A_{22}+A_{23}, 
\end{align*}
where we define  
\begin{align*}
A_{21}= &~ \big|P\left(\Theta\in I_{\delta}\mid \mathcal{E}_{n,k}\right)-P\left(\Theta\in I_{\delta}\cap I_{k,\eta}\mid \mathcal{E}_{n,k} \right) \big|=P\left(\Theta\in I_{\delta}\cap I_{k,\eta}^c\mid \mathcal{E}_{n,k} \right), \notag\\
A_{22}= &~\int |P_{n,i}(\theta_k)-P_{n,i}(\theta)|\, f_{n,k,-i}(\theta)\, I(\theta\in I_{\delta}\cap I_{k,\eta})\, \md \theta, \notag\\
A_{23} = &~ P(Y_{n,i}=1, \Theta\in I_{\delta} \cap I_{k,\eta}^c\mid \mathcal{E}_{n,k} ),
\end{align*}
where in the definition of $A_{22}$, we let  $f_{n,k,-i}(\theta)$ denote the probability density of $\Theta$ conditioning on $\mathcal{E}_{n,k}$. 
Since $A_{23}\leqslant A_{21}
$ and $A_2$ and $A_2^*$ can be analyzed similarly, to prove \eqref{eq:a2a2starbd}, it suffices to show  $A_{21}\leqslant \epsilon/12$ and $A_{22}\leqslant \epsilon/12$ below. 
% derive upper bounds of $A_{21}$ and $A_{22}$ below.  

First, $A_{21}=A_{21,num}/P( \mathcal{E}_{n,k})$, where we define $A_{21,num}=P(\Theta\in I_{\delta}\cap I_{k,\eta}^c, \mathcal{E}_{n,k})$ satisfying
\begin{align}
&~A_{21,num}\notag\\
=& P\left( |\Theta-\theta_k|>n^{-\eta},\  \Theta \in I_{\delta},\ \bar{Y}_{n,-i}=\frac{k}{n-1}\right)\notag\\
=& P\left( \big|\Theta-\bar{P}_{n,-i}^{-1}\big(\bar{Y}_{n,-i}\big)\big|>n^{-\eta},\  \Theta \in I_{\delta},\ \bar{Y}_{n,-i}=\frac{k}{n-1}\right)\hspace{1.5em}(\text{by }\bar{P}_{n,-i}(\theta_k)=k/(n-1))\notag\\
%= & \int P\left( \big|\theta-\bar{P}_{n,-i}^{-1}\big(\bar{Y}_{n,-i}\big)\big|> n^{-\alpha},\ \bar{Y}_{n,-i}=\frac{k}{n-1} \mid \Theta=\theta\right) \indicator( \theta\in I_{\delta} ) \md \theta \notag\\
= & \int P\left\{ \bar{P}_{n,-i}^{-1}(\bar{Y}_{n,-i})>\theta+n^{-\eta} \text{ or }<\theta-n^{-\eta},\ \bar{Y}_{n,-i}=\frac{k}{n-1} \mid \Theta=\theta\right\} \indicator( \theta\in I_{\delta} ) \md  \theta \notag\\
\leqslant & \int P\left\{ \bar{Y}_{n,-i}>\bar{P}_{n,-i}(\theta+n^{-\eta}) \text{ or }<\bar{P}_{n,-i}(\theta-n^{-\eta}),\ \bar{Y}_{n,-i}=\frac{k}{n-1} \ \biggr| \ \Theta=\theta\right\} \indicator( \theta\in I_{\delta} ) \md \theta \notag\\
\leqslant & \int P\left\{ \big|\bar{Y}_{n,-i}-\bar{P}_{n,-i}(\theta)\big|>|\bar{P}_{n,-i}'(\theta_{\eta})| n^{-\eta} \mid \Theta=\theta\right\} \indicator( \theta\in I_{\delta} ) \md \theta \label{eq:i1numerator}
%\leqslant & \int P\left\{ \bar{Y}_{n,-i}-P_{n,-i}(\theta)>P_{n,-i}'(\theta_{\alpha}) n^{-\alpha} \text{ or }<-P_{n,-i}'(\tilde{\theta}_{\alpha}) n^{-\alpha} \mid \Theta=\theta\right\} \indicator( \theta\in I_{\delta} ) \md \theta \label{eq:i1numerator}
\end{align}
%where in the second equation,
%we use $\bar{P}_{n,-i}(\theta_k)=\frac{k}{n-1}$,
where  the last inequality is obtained by the intermediate value theorem,   and 
 $\theta_{\eta}$ is between 
$\theta-n^{-\eta}$ and $\theta+n^{-\eta}$. 
%$\theta$ and $\theta+n^{-\alpha}$, and $\tilde{\theta}_{\alpha}$ is between $\theta$ and $\theta-n^{-\alpha}$.  
When $n$ is sufficiently large, 
 $\theta_{\eta}\in (\delta/2, 1-\delta/2)$  given $\theta\in(\delta,1-\delta)$. 
%given $\theta\in(\delta,1-\delta)$, 
%we know $\theta_{\alpha}\in (\delta/2, 1-\delta/2)$ 
Thus, by Condition \ref{cond:derivbound},  
there exists a constant $m_{\delta/2}>0$  independent with $(n,i)$ such that
$|P_{n,-i}'(\theta_{\eta})|>m_{\delta/2}$. 
Therefore, 
\begin{align*}
\eqref{eq:i1numerator} \leqslant\int P\left\{ \big|\bar{Y}_{n,-i}-\bar{P}_{n,-i}(\theta)\big|>m_{\delta/2} n^{-\eta} \mid \Theta=\theta\right\} \indicator( \theta\in I_{\delta} ) \md \theta.   	
\end{align*}
By  Lemma \ref{lm:concenybarexp}, i.e., Hoeffding's inequality of bounded variables, we have
%the Hoeffding's inequality of bounded variables (see Lemma \ref{lm:concenybarexp}), we have
%concentration inequality in Lemma \ref{lm:concenybarexp}, we have
\begin{align*}
	\eqref{eq:i1numerator} \leqslant 2\int e^{- 2(n-1)n^{ -2\eta} m_{\delta/2}^2 } \, \indicator(\theta\in I_{\delta}) \md \theta \leqslant 2\delta e^{-  2(n-1)n^{ -2\eta} m_{\delta/2}^2 }. 
\end{align*}
%we can take $\theta_{\alpha}\to 0$ such that 
%$n^{1-2\alpha} (P_n'(\theta_{\alpha}))^2\to \infty $, so the above bound $\to 0$ after timed by the denominator bound.  
Combining the above inequality with   Lemma \ref{lm:lowbdofbary}, we have
\begin{align*}
	A_{21}=\frac{A_{21,num}}{P( \mathcal{E}_{n,k})}\leqslant 2\delta e^{- (n-1)n^{ -2\eta} m_{\delta/2}^2 }\, \frac{n}{ \tilde{C}_{\alpha\beta}}. 
\end{align*} 
When $0<\eta<1/2$, the exponential term converges to 0 faster  than $n^{-r}$ for any   $r>0$. Thus, for any $\epsilon>0$, there exists $N_{\epsilon,\alpha,\beta,\delta}$ such that for $n\geqslant N_{\epsilon,\alpha,\beta,\delta}$, $A_{21}\leqslant \epsilon/12$.

Second, by  $\theta_k\in (\alpha,\beta)$, when $n$ is sufficiently large, 
$I_{k,\eta}\subseteq (\alpha/2,(1+\beta)/2)$.
Therefore,  
\begin{align*}
	\sup_{\theta\in I_{k,\eta}}\left|P_{n, i}(\theta)-P_{n, i}\left(\theta_k\right)\right|\leqslant \sup_{\eta\in I_{k,\eta}}|P'_{n,i}(\eta)| \sup_{\theta\in I_{k,\eta}}|\theta-\theta_k|\leqslant M_{\alpha\beta,2}n^{-\eta},
\end{align*} 
where $ M_{\alpha\beta,2}$ is a constant that is independent of $(n,i)$ by Condition \ref{cond:derivbound}. 
It follows that
\begin{align*}
	A_{22}\leqslant &~ \int M_{\alpha\beta,2}n^{-\eta} f_{n, k,-i}\left(\theta\right)\indicator(\theta\in I_{\delta}\cap I_{k,\eta}) \md \theta\notag\\ 
 \leqslant &~ M_{\alpha\beta,2}n^{-\eta} \times P(\Theta\in I_{\delta}\cap I_{k,\eta} \mid \mathcal{E}_{n,k} )\notag\\ 
 \leqslant &~ M_{\alpha\beta,2} n^{-\eta} \leqslant \epsilon/12, 
\end{align*}
 when $n$ is sufficiently large. 

\begin{remark}\label{rm:novelty}
Although Theorem \ref{thm:asymepsilonbd}     concerns a strict subset of $(0,1)$, 
its proof   markedly differs from that in \cite{douglas2001asymptotic}, since  Assumptions 4--5 in \cite{douglas2001asymptotic} cannot be directly applied. 
Notably, 
  in  Step 2 of the preceding proof, an interval $I_{\delta}$ is introduced to establish the necessary  inequalities, which is not required in \cite{douglas2001asymptotic}. 
  In addition, Proposition 1 in \cite{douglas2001asymptotic}  shows that an IRF sequence adhering to their Assumptions 1--5 can be approximated by another IRF sequence satisfying Conditions \ref{cond:unilocal}--\ref{cond:unifbd} above. 
 However, this approximation result does not immediately imply that the nonparametric IRF sequence can be uniquely identified  within the extended model class.  
 \end{remark}

\section{Discussions}\label{sec:disc}

Nonparametric IRT models provide a versatile   framework and play an important role in  ensuring robust  measurement  against model misspecification and in assessing parametric IRT model fit. 
In practical applications, it is often imperative to consider a large space of functions that  embraces popular parametric IRT models.
This requirement, however, considerably   complicates the study of model identifiability.  
In this paper, we show that 
 the assumptions  restricting the model class in \cite{douglas2001asymptotic} can be substantially relaxed, and we 
 establish  asymptotic identifiability  for an extended model class that includes many popular parametric IRT models.   
  The results serve as a crucial link between parametric and nonparametric IRT models and  provide a solid theoretical foundation for using nonparametric IRT models in  large-scale assessments with many items. 

%  nicely bridge the parametric and nonparametric IRT models and provide  solid theoretical justification for the applications of  nonparametric IRT models for  assessments with many items.    

\section*{Acknowledgements}
The research is partially supported by Wisconsin Alumni Research Foundation. 
The author would like to thank Prof. Zhiliang Ying  for helpful discussions.

\section*{Appendix}

We present all the lemmas in Section A  and provide their proofs in Section B.
The proofs of propositions are provided in Section C.  
 
\subsection*{A. Lemmas}
\begin{lemma}\label{lm:lowbdofbary}
Consider an integer $k$ such that $\theta_k\in (\alpha,\beta)$. 
There exist   constants $\tilde{C}_{\alpha\beta }>0$ and $N_{\alpha\beta}$ independent with $(n,i)$ 
such that when $n \geqslant N_{\alpha\beta}$, 
%where $N_{\alpha\beta}$ is an integer independent with $i$, such that  
%is sufficiently large, 
\begin{align*} 
	P\left(\bar{Y}_{n,-i}=\frac{k}{n-1}\right)\geqslant \tilde{C}_{\alpha\beta}n^{-1}.
\end{align*}	
\end{lemma}

\begin{lemma}\label{lm:fixedrange}
Under the conditions of Theorem \ref{thm:asymepsilonbd}, 
%when $\delta$ is sufficiently small, 
there exist  constants  $\tilde{C}_{\alpha \beta,1}$  and $\tilde{C}_{\alpha \beta,2}$  independent with $(n,i)$  such that 
% $C_{\alpha,\beta,2}>0$ and $\delta_{\alpha,\beta}>0$ independent with $(n,i)$ such that when $\delta\leqslant \delta_{\alpha,\beta}$, 
\begin{align*}
	\left|1-P\left(\Theta \in I_{\delta}\mid \mathcal{E}_{n,k} \right)\right|\leqslant \frac{ n\delta\exp(-n\tilde{C}_{\alpha \beta,1} )}{\tilde{C}_{\alpha\beta,2}}. 
\end{align*}
%where $\tilde{C}_{\alpha\beta,1}$ is specified in Lemma \ref{lm:lowbdofbary}. 
%\begin{align} %\label{eq:fixedrangeapprx}
%	\left|1-P\left(\Theta \in I_{\delta}\mid \mathcal{E}_{n,k} \right)\right|\leqslant C_{\alpha,\beta,2}\, {\max_{j\neq i} P_{n,j}(\delta) \delta} . \notag
%\end{align}  	
\end{lemma}

\begin{lemma}[Hoeffding inequality of bounded variables]
\label{lm:concenybarexp}
For any $(n,i)$, and $m>0$, 
	\begin{align*}
		P(|\bar{Y}_{n,-i}-\bar{P}_{n,-i}(\theta)|> m \mid \Theta=\theta)\leqslant 2\exp[- 2(n-1)m^2]. 
	\end{align*}
\end{lemma} 
%\begin{proof}
%Lemma \ref{lm:concenybarexp} follows from \cite{hoeffding1994probability}, where we note that $\bar{Y}_{n,-i}$ is an average over $n-1$ Bernoulli random variables, and $2(n-1) \geqslant  n$. 	
%\end{proof}

%\section*{Appendix. Proofs of Lemmas}
\vspace{1em}
\subsection*{B. Proofs of Lemmas} 
\subsubsection*{B1. Proof of Lemma \ref{lm:lowbdofbary}}
%Let $\mu_\theta$ and $\sigma_\theta^2$ denote the mean and variance of $(n-1)\bar{Y}_{n,-i}$
%%$n \bar{Y}_n=\sum_{j=1}^n Y_{n,j}$ 
%conditional on $\theta$.
Let $\mu_{\theta}=\sum_{j\neq i} P_{n,j}(\theta)$ and $\sigma^2_{\theta}=\sum_{j\neq i} P_{n,j}(\theta)(1-P_{n,j}(\theta))$, 
which are the mean and variance of $(n-1)\bar{Y}_{n,-i}$ conditioning on $\Theta=\theta$, respectively.  
%We point out that $\mu_{\theta}$  
By applying a bound
on the normal approximation for the distribution of a sum of independent Bernoulli variables
\citep{mikhailov1994refinement} 
\begin{align}
	\left|P\left[(n-1)\bar{Y}_{n,-i}=k  \mid \Theta=\theta\right]-\frac{1}{\sigma_\theta \sqrt{2 \pi}} e^{-\frac{\left(k-\mu_\theta+1 / 2\right)^2}{2 \sigma_\theta^2}}\right| \leq \frac{c}{\sigma_\theta^2}\label{eq:normalappx}
\end{align}
 for some universal constant $c$. 

Define the intervals $I_{n,1/2}=(\theta_k-n^{-1/2}, \theta_k+n^{-1/2})$ and $\tilde{I}_{n,1/2}=I_{n,1/2}\cap (\alpha/2,(1+\beta)/2)$. 
For $\theta \in \tilde{I}_{n,1/2}$ and $\theta_k \in (\alpha,\beta)$, 
%Define the interval $I_{n,1/2}=(\theta_k-n^{-1/2}, \theta_k+n^{-1/2})$. 
%Given $\theta_k\in (\alpha,\beta)$,
%$I_{n,1/2}\subseteq (a/2,(1+b)/2)$ when  $n$ is sufficiently large. 
%When $n$ is sufficiently large, $I_{n,1/2}\subseteq $
%Given $\theta_k\in (\alpha,\beta)$, we take  $\theta\in(\theta_k-n^{-1/2}, \theta_k+n^{-1/2}) \subseteq (a/2,(1+b)/2)$ when $n$ is sufficiently.
% such that $|\theta-\theta_k|<n^{-1/2}$. 
%There exists constants $\tilde{a}$ and $\tilde{b}$ such that 
%When $n$ is sufficiently large, $\theta\in (a/2,(1+b)/2)$. 
%$\theta\in (\tau/2,1-\tau/2)$. 
%Thus for $\theta\in I_{n,1/2}$, 
\begin{align}
	\frac{|k-\mu_{\theta}|}{n-1}=|\bar{P}_{n,-i}(\theta_k)-\bar{P}_{n,-i}(\theta)|\leqslant {M}_{\alpha\beta,2} |\theta_k-\theta| < {M}_{\alpha\beta,2}n^{-1/2}.\label{eq:kmuthetadiff}
\end{align} 
where in the second inequality,  ${M}_{\alpha\beta,2}$ is  a constant that is independent with $(n,i)$  by Condition \ref{cond:derivbound}. 
Moreover, by Condition \ref{cond:unifbd}, 
%\blue{(assumption?)},  
there exist positive constants $L_{\alpha\beta} $ and  $U_{\alpha\beta} $  that are independent with $(n,i)$ such that 
\begin{align}
L_{\alpha\beta} <\sigma^2_{\theta}/(n-1) < U_{\alpha\beta}.\label{eq:sigmathetabd}
\end{align}
%Thus $\theta$ is also bounded away from 0 and 1 when $n$ is sufficiently large.  
Therefore, 
%when  $n$ is sufficiently large, 
by \eqref{eq:normalappx}, \eqref{eq:kmuthetadiff}, and \eqref{eq:sigmathetabd}, for $\theta\in \tilde{I}_{n,1/2}$, 
\begin{align} %\label{eq:denominatorlb}
	P\left(\bar{Y}_{n,-i}=\frac{k}{n-1} \ \Big| \   \theta\right) > &~ \frac{1}{\sqrt{2\pi\sigma^2_{\theta} }} e^{-\frac{\left(k-\mu_\theta+1 / 2\right)^2}{2\sigma_\theta^2}} - \frac{c}{\sigma^2_{\theta}} \notag\\
>  &~ \frac{1}{\sqrt{2U_{\alpha\beta}(n-1)\pi}} e^{-\frac{\left({M}_{\alpha\beta,2}(n-1)n^{-1/2} + \frac{1}{2}\right)^2}{2 L_{\alpha\beta}(n-1)}} - \frac{c}{L_{\alpha\beta}(n-1)}\notag\\
 > &~ \frac{\tilde{C}_{\alpha\beta}}{2n^{1/2}}, \notag 
\end{align}
where $\tilde{C}_{\alpha\beta}>0$ is a constant that is independent with $(n,i)$. 
Therefore, 
\begin{align} %\label{eq:denominatorlb}
P\left(\bar{Y}_{n,-i}=\frac{k}{n-1} \right) \geqslant &~\int_{\theta\in \tilde{I}_{n,1/2}}	P\left(\bar{Y}_{n,-i}=\frac{k}{n-1}  \ \Big| \   \theta\right) \md\theta >  \int_{\theta \in \tilde{I}_{n,1/2}}\frac{\tilde{C}_{\alpha\beta}}{2n^{1/2}}\md\theta.\label{eq:lbprobdeno}
\end{align}
There exists $N_{\alpha\beta}$ independent with $i$ such that when $n\geqslant N_{\alpha\beta}$,
%When $n$ is sufficiently large, 
   $\tilde{I}_{n,1/2}=I_{n,1/2}$ by $\theta_k\in (\alpha,\beta)$. 
Then by \eqref{eq:lbprobdeno}, 
\begin{align*}
	P\left(\bar{Y}_{n,-i}=\frac{k}{n-1} \right)>\int_{\theta \in  {I}_{n,1/2}}\frac{\tilde{C}_{\alpha\beta}}{2n^{1/2}}\md\theta = \tilde{C}_{\alpha\beta}n^{-1}. 
\end{align*}

\vspace{1em}
\subsubsection*{B2. Proof of  Lemma  \ref{lm:fixedrange}}
To prove Lemma \ref{lm:fixedrange}, we note that
\begin{align*}
1-P(\Theta \in I_{\delta}\mid \mathcal{E}_{n,k})={P\left(0<\Theta<\delta\mid  \mathcal{E}_{n,k}\right)+P\left(1-\delta <\Theta<1 \mid  \mathcal{E}_{n,k}\right)}. 
\end{align*} 
We next derive an upper bound of $P\left(0<\Theta<\delta\mid  \mathcal{E}_{n,k}\right)$,
 and a similar upper bound can be obtained for $P\left(1-\delta <\Theta<1 \mid  \mathcal{E}_{n,k}\right)$ following a similar analysis. 
  
By $P(0<\Theta<\delta\mid  \mathcal{E}_{n,k})= P(0<\Theta<\delta,\ \mathcal{E}_{n,k})/P(\mathcal{E}_{n,k})$, 
%We next first derive an upper bound of $P(0<\Theta<\delta,\ \mathcal{E}_{n,k})$.
and the lower bound of $P(\mathcal{E}_{n,k})$ in Lemma \ref{lm:lowbdofbary},
it suffices to derive an upper bound of $P(0<\Theta<\delta,\ \mathcal{E}_{n,k})$ below. In particular,  
\begin{align}
	 P\left(0<\Theta<\delta,\ \mathcal{E}_{n,k}\right) 
	=\int P\biggr(\bar{Y}_{n,-i}=\frac{k}{n-1}\mid \Theta=\theta \biggr) \, \indicator(0<\theta<\delta)\mathrm{d}\theta. \label{eq:jointprobintegral}  
\end{align}
%Let $\bar{f}_{-i}=\sum_{j\neq i} f_j/(n-1)$. 
By $\bar{P}_{n,-i}(\theta_k)=k/(n-1)$, 
\begin{align}
\bar{Y}_{n,-i}=\frac{k}{n-1}\quad \Leftrightarrow \quad \bar{Y}_{n,-i}-\bar{P}_{n,-i}(\theta)=[\bar{P}_{n,-i}(\theta_k)-\bar{\kappa}_{-i}]-[\bar{P}_{n,-i}(\theta)-\bar{\kappa}_{-i}],\label{eq:equicenterevent}
\end{align}
where we define $\bar{\kappa}_{-i}=\sum_{j\neq i} \kappa_j/(n-1)$. 
%For $\theta\in [0,\delta]$, 
%as IRFs are increasing functions by Condition \ref{cond:derivbound}, 
%\begin{align}\label{eq:IRFsdeltabd}
%	P_{n,-i}(\theta)\leqslant P_{n,-i}(\delta). 
%\end{align}
%\begin{align*}
%P\biggr(\bar{Y}_{n,-i}=\frac{k}{n-1}\mid \Theta=\theta \biggr) = P\biggr(\bar{Y}_{n,-i}-\bar{P}_{n,-i}(\theta)=\frac{k}{n-1}-\bar{f}_{-i}+\bar{f}_{-i}-\bar{P}_{n,-i}(\theta)\mid \Theta=\theta \biggr) \notag 
%\end{align*}
%\begin{align*}
%&~P\biggr(\bar{Y}_{n,-i}=\frac{k}{n-1}\mid \Theta=\theta \biggr) \notag\\
%	\leqslant &~ P\biggr(\bar{Y}_{n,-i}-\bar{P}_{n,-i}(\theta)=\frac{k}{n-1}-f+f-\bar{P}_{n,-i}(\theta)\mid \Theta=\theta \biggr) \notag 
%\end{align*}
%By $\bar{P}_{n,-i}(\theta_k)=k/(n-1)$, 
%\begin{align*}
%	\frac{k}{n-1}-\bar{f}_{-i} = \bar{P}_{n,-i}(\theta_k)-\bar{f}_{-i} \geqslant \bar{P}_{n,-i}(a)-\bar{P}_{n,-i}(a/2) 
%\end{align*}
%where in the last inequality, we use 
%$\theta_k\in (\alpha,\beta)$ and Condition \ref{cond:unifbd}
By $0<\alpha/2<\alpha<\theta_k<\beta$ and Conditions \ref{cond:derivbound} and  \ref{cond:unifbd}, we have $\bar{P}_{n,-i}(\theta_k) \geqslant \bar{P}_{n,-i}(\alpha)$ and $\bar{P}_{n,-i}(\alpha/2)\geqslant \bar{\kappa}_{-i}.$  
%\begin{align}\label{eq:lwbdaf}
%\bar{P}_{n,-i}(\theta_k) \geqslant \bar{P}_{n,-i}(a)\hspace{1em}\text{and}\hspace{1em} 	\bar{P}_{n,-i}(a/2)\geqslant \bar{f}_{-i}. 
%\end{align}
Therefore, 
\begin{align} 
\bar{P}_{n,-i}(\theta_k)-\bar{\kappa}_{-i} \geqslant 	&~\bar{P}_{n,-i}(\alpha)-\bar{P}_{n,-i}(\alpha/2) = \bar{P}_{n,-i}'(\tilde{\alpha} )\alpha/2\geqslant \tilde{m}_{\alpha} \alpha/2, \label{eq:knratiolowerbd} %\notag\\
%=&~\bar{P}_{n,-i}'(\tilde{a} )a/2 \geqslant \tilde{m}_a a/2, \label{eq:knratiolowerbd}
\end{align}
where the second equation is obtained by the intermediate value theorem with  $\tilde{\alpha}\in (\alpha/2,\alpha)$,
and the last inequality is obtained by Condition \ref{cond:derivbound} with $\tilde{m}_{\alpha}$ being a constant independent with $(n,i)$.
Let $\tilde{C}_{\alpha} = \tilde{m}_{\alpha} \alpha/2$.  
By Condition \ref{cond:unifbd}, 
there exists $\delta_{\alpha}>0$ such that for  $\theta<\delta\leqslant \delta_{\alpha}$, 
\begin{align}
	 \bar{P}_{n,-i}(\theta) -\bar{\kappa}_{-i} < \bar{P}_{n,-i}(\delta) -\bar{\kappa}_{-i} < \tilde{C}_{\alpha} /2.  \label{eq:IRFsdeltabd}
\end{align}
Combining \eqref{eq:equicenterevent}, \eqref{eq:knratiolowerbd}, and \eqref{eq:IRFsdeltabd},  
\begin{align}
&~P\left( \bar{Y}_{n,-i}=\frac{k}{n-1} \mid \Theta = \theta\right)\notag\\
 = &~ P\left( \bar{Y}_{n,-i}-\bar{P}_{n,-i}(\theta)=[\bar{P}_{n,-i}(\theta_k)-\bar{\kappa}_{-i}]-[\bar{P}_{n,-i}(\theta)-\bar{\kappa}_{-i}]\mid \Theta = \theta\right)\notag\\
\leqslant &~P\left( \bar{Y}_{n,-i}-\bar{P}_{n,-i}(\theta)\geqslant  \tilde{C}_{\alpha} /2 \mid \Theta=\theta\right)	\notag\\
\leqslant &~ 2\exp[-(n-1)\tilde{C}_{\alpha}^2 /2  ], \notag 
%\leqslant &~\frac{ 4}{\tilde{C}_a^2} E(\bar{Y}_{n,-i}-\bar{P}_{n,-i}(\theta))^2\hspace{3em} (\text{by Markov's inequality})\notag\\
%=&~\frac{4}{\tilde{C}_{\alpha}^2(n-1)^2}\sum_{j\neq i}P_{n,j}(\theta)(1-P_{n,j}(\theta))\notag\\
%\leqslant &~ \frac{8 \max_{j\neq i}P_{n,j}(\delta) }{n\tilde{C}_{\alpha}^2}. \notag
\end{align}  
where the last inequality follows by Lemma \ref{lm:concenybarexp}. 
By the above inequality and 
  \eqref{eq:jointprobintegral},
 \begin{align*}
		P\left(0<\Theta<\delta \mid  \mathcal{E}_{n,k}\right)\leqslant\frac{2\delta\exp[-(n-1)\tilde{C}_{\alpha}^2 /2  ]}{ P(\mathcal{E}_{n,k}) }\leqslant \frac{2n\delta\exp[-(n-1)\tilde{C}_{\alpha}^2 /2  ]}{\tilde{C}_{\alpha\beta}},  
%		\\=\frac{8\max_{j\neq i}P_{n,j}(\delta)\delta}{n\tilde{C}_{\alpha}^2}\frac{1}{P(\mathcal{E}_{n,k})}\leqslant \frac{8\max_{j\neq i} P_{n,j}(\delta) \delta}{\tilde{C}_{\alpha\beta} \tilde{C}_{\alpha}^2}, 
\end{align*}
%\begin{align*}
%		P\left(0<\Theta<\delta \mid  \mathcal{E}_{n,k}\right)\leqslant \frac{8\max_{j\neq i}P_{n,j}(\delta)\delta}{n\tilde{C}_{\alpha}^2}\frac{1}{P(\mathcal{E}_{n,k})}\leqslant \frac{8\max_{j\neq i} P_{n,j}(\delta) \delta}{\tilde{C}_{\alpha\beta} \tilde{C}_{\alpha}^2}, 
%\end{align*}
where the second inequality is obtained by Lemma \ref{lm:lowbdofbary}. 
A similar upper bound can be obtain for   
$P (1-\delta <\Theta<1 \mid  \mathcal{E}_{n,k}\ ) $  too. 
Lemma \ref{lm:fixedrange} is proved. 

%Combining with Lemma \ref{lm:lowbdofbary}, we obtain
%\begin{align*}
%	P\left(0<\Theta<\delta \mid  \mathcal{E}_{n,k}\right)\leqslant \frac{8\max_{j\neq i} P_{n,j}(\delta) \delta}{n \tilde{C}_{a}^2 \tilde{C}_{\alpha\beta}n^{-1}} = \frac{8\max_{j\neq i} P_{n,j}(\delta) \delta}{\tilde{C}_{\alpha\beta} \tilde{C}_{a}^2}. 
%\end{align*} 

%Following similar analysis, we can establish a similar upper bound of $ P(1-\delta<\Theta\leqslant 1 \mid P(\mathcal{E}_{n,k}) $. 
%Lemma \ref{lm:fixedrange} follows by 

\vspace{1em}
\subsection*{C. Proofs of Propositions} 

\subsubsection*{C1. Proof of Proposition \ref{prop:normalogive}} Let $x=\Phi^{-1}(\theta)$. We can equivalently write $P'_{n,i}(\theta)=h_{n,i}(x)$, where 
\begin{align*}
h_{n,i}(x)=\frac{a_{n,i} \phi[a_{n,i}(x-b_{n,i})] }{ \phi(x)} = a_{n,i}\exp\left[-\frac{1}{2}(a_{n,i}^2-1)x^2+a_{n,i}^2b_{n,i}\left(x-\frac{1}{2}\right) \right]. 
\end{align*}
Note that $\theta \downarrow 0$ and  $\theta \uparrow 0$ correspond to $x \to -\infty$ and $x\to +\infty$, respectively. 
Let $p_{+}=\lim_{\theta\to 1}P_{n,i}'(\theta)=\lim_{x\to +\infty} h_{n,i}(x)$ and $p_{-}=\lim_{\theta\to 0}P_{n,i}'(\theta)=\lim_{x\to -\infty} h_{n,i}(x)$.  
As $h_{n,i}(x)=0$ when $a_{n,i}=0$,
it suffices to consider $a_{n,i}\neq 0$ below. In particular, 
\begin{align*}
(p_{-}, p_{+})=\begin{cases}
(+\infty, +\infty) & \text{ when } 0<a_{n,i}^2<1 \\
(0, 0) & \text{ when } a_{n,i}^2>1\\
(0, +\infty) & \text{ when } a_{n,i}^2=1, b_{n,i}>0\\
(  +\infty, 0) & \text{ when } a_{n,i}^2=1, b_{n,i}<0\\
(1, 1) & \text{ when } a_{n,i}^2=1, b_{n,i}=0. 
\end{cases}	
\end{align*}
This suggests that $P'_{n,i}(\theta)$ cannot be uniformly bounded on $(0,1)$ when $a_{n,i}^2\neq 1$ or $ b_{n,i}\neq 0. $

On the other hand, when  
$\theta\in [\alpha, \beta]$, $x\in [\Phi^{-1}(\alpha), \Phi^{-1}(\beta)]$ with the two end points bounded away from $-\infty$ and $+\infty$.  
Therefore, when $a_{n,i}\in  [1/C, C]$ and $|b_{n,i}|\in [1/C', C']$,
there exist  $0< m_{\alpha \beta} < M_{\alpha \beta} <\infty$ such that 
$P'_{n,i}(\theta) = h_{n,i}(x) \in (m_{\alpha \beta}, M_{\alpha \beta})$.  
Thus Condition \ref{cond:derivbound} is satisfied. 

We next prove that Condition \ref{cond:unifbd} is also satisfied. 
When $\epsilon\in (0,1/2)$, $\Phi^{-1}(\epsilon)<0$, and then 
 we set $l_{\epsilon}=\Phi\left[C_a\Phi^{-1}(\epsilon) - C_b\right].$ 
When $\theta\in[0,l_{\epsilon}]$, 
%\begin{align*}
%	l_{\epsilon}=\Phi\left[C\Phi^{-1}(\epsilon) + C'\right]. 
%\end{align*}
\begin{align*}
	P_{n,i}(\theta)\leqslant &~ P_{ni}(l_{\epsilon}) =\Phi\big[a_{ni}(\Phi^{-1}(l_{\epsilon}) - b_{ni})\big] 
	\leqslant   \Phi\{ a_{ni}[C_a\Phi^{-1}(\epsilon)-C_b+ |b_{ni}| ] \}\notag\\
	\leqslant &~ \Phi\{ a_{ni}[C_a\Phi^{-1}(\epsilon)-C_b+C_b] \} \leqslant \Phi\left[\frac{1}{C_a} C_a \Phi^{-1}(\epsilon)\right]= \epsilon. 
\end{align*}
%When $\epsilon\in (0,1/2)$, $\Phi^{-1}(\epsilon)<0$ and then, 
%\begin{align*}
%	 a_{ni}\big[C_a\Phi^{-1}(\epsilon)-C_b-b_{ni} \big] \leqslant &~a_{ni} \big[C_a\Phi^{-1}(\epsilon)-C_b+|b_{ni}| \big]  
%	 \leqslant a_{ni}C_a\Phi^{-1}(\epsilon) \leqslant \Phi^{-1}(\epsilon). 
%\end{align*}
When $\epsilon \in (1/2,1)$, $\Phi^{-1}(\epsilon)>0$, and  
 we set $l_{\epsilon}=\Phi\left[C_a^{-1}\Phi^{-1}(\epsilon) - C_b\right].$
Then  we  can  obtain $P_{n,i}(\theta)\leqslant \epsilon$ similarly to the above analysis. 
Following similar analysis, we can construct $u_{\epsilon}$ so that $1-P_{n,i}(\theta)\leqslant \epsilon$ for $\theta \in [u_{\epsilon},1]$. %An upper bound $u_{\epsilon}$ can be constructed similarly. 
In brief, Condition \ref{cond:unifbd} is satisfied. 

\subsubsection*{C2. Proof of Proposition \ref{prop:4pllogisticmodel}}
Let $x=\Phi^{-1}(\theta)$. We can equivalently write $P'_{n,i}(\theta)=h_{n,i}(x)$, where 
\begin{align*}
h_{n,i}(x)=&~(d_{n,i}-c_{n,i})a_{n,i}\frac{ g'\left[a_{n,i}(x-b_{n,i})\right]}{\phi(x)} \\
=&~ (d_{n,i}-c_{n,i})a_{n,i}  \frac{ e^{x^2/2}/[e^{a_{n,i}(x-b_{n,i})}+e^{-a_{n,i}(x-b_{n,i})}]}{1+2/[e^{a_{n,i}(x-b_{n,i})}+e^{-a_{n,i}(x-b_{n,i})}]}   
%=&~ (d_{n,i}-c_{n,i})a_{n,i} \frac{1}{2}\left[1-\frac{1}{1+2/[e^{a_{n,i}(x-b_{n,i})}+e^{-a_{n,i}(x-b_{n,i}]})} \right]e^{x^2/2}  
\end{align*}
which is obtained by %$g'(x)=\frac{1}{e^x+e^{-x}+2}=\frac{1}{2}()$
\begin{align*}
	g'(x)=  \frac{1/(e^x+e^{-x})}{1+2/(e^x+e^{-x})}. 
\end{align*} 
Note that $\theta \downarrow 0$ and  $\theta \uparrow 0$ correspond to $x \to -\infty$ and $x\to +\infty$, respectively. 
%Let $p_{+}=\lim_{\theta\uparrow 1}P_{n,i}'(\theta)=\lim_{x\to +\infty} h_{n,i}(x)$ and $p_{-}=\lim_{\theta\downarrow 0}P_{n,i}'(\theta)=\lim_{x\to -\infty} h_{n,i}(x)$.  
When $a_{n,i}=0$, $h_{n,i}(x)=0$.
When $a_{n,i}\neq 0$,   $e^{a_{n,i}(x-b_{n,i})}+e^{-a_{n,i}(x-b_{n,i})}\to +\infty$ and $e^{x^2/2}$ as $|x|\to +\infty$. 
Moreover, as the quadratic term diverges faster than the linear term, we know $e^{x^2/2}/[e^{a_{n,i}(x-b_{n,i})}+e^{-a_{n,i}(x-b_{n,i})}]\to \infty$. 
Thus,   $h_{n,i}(x)\to +\infty$.
 
On the other hand, when  
$\theta\in [\alpha, \beta]$, $x\in [\Phi^{-1}(\alpha), \Phi^{-1}(\beta)]$ with the two end points bounded away from $-\infty$ and $+\infty$.  
Therefore, under the conditions in (ii) of Proposition \ref{prop:4pllogisticmodel}, 
 there exist  $0< m_{\alpha \beta} < M_{\alpha \beta} <\infty$ such that 
$P'_{n,i}(\theta) = h_{n,i}(x) \in (m_{\alpha \beta}, M_{\alpha \beta})$.  
Thus, Condition \ref{cond:derivbound} is satisfied. 

We next prove that Condition \ref{cond:unifbd} is also satisfied. 
When $\epsilon$ is small such that  $g^{-1}(\epsilon/C_{c,d})<0$, we set $l_{\epsilon}=\Phi\left[C_ag^{-1}(\epsilon/C_{c,d}) - C_b\right].$ Then for $\theta\in [0,l_{\epsilon}]$, 
\begin{align*}
%	P_{n,i}(\theta)-c_{n,i}\leqslant &~\left(d_{n,i}-c_{n,i}\right) g\left[a_{n,i}\left(\Phi^{-1}(l_{\epsilon})-b_{n,i}\right)\right] \\
	P_{n,i}(\theta)-c_{n,i} \leqslant P_{n,i}(l_{\epsilon})-c_{n,i} = &~ \left(d_{n,i}-c_{n,i}\right) g\left(a_{n,i}[C_ag^{-1}(\epsilon/C_{c,d}) - C_b - b_{n,i} ] \right) \\
	\leqslant &~   \left(d_{n,i}-c_{n,i}\right) g\left(a_{n,i}[C_ag^{-1}(\epsilon/C_{c,d}) - C_b + |b_{n,i}|] \right) \\
 \overset{(i1)}{\leqslant} &~ \left(d_{n,i}-c_{n,i}\right) g[g^{-1}(\epsilon/C_{c,d})]   
	\leqslant   \epsilon, 
\end{align*} 
where    the inequality $(i1)$ above is obtained by $a_{n,i}C_a\geqslant 1 > 0$ and $|b_{n,i}|-C_b\leqslant 0$. 
When $\epsilon$ is large such that $g^{-1}(\epsilon/C_{c,d})>0$,
we set $l_{\epsilon}=\Phi\left[C_a^{-1}g^{-1}(\epsilon/C_{c,d}) - C_b\right].$ and then $P_{n,i}(\theta)-c_{n,i}\leqslant \epsilon$ can be obtained similarly. 
Following similar analysis, we can construct $u_{\epsilon}$ so that 
$d_{n,i}-P_{n,i}(\theta) \leqslant \epsilon $ for $\theta \in [u_{\epsilon},1]$. In brief, Condition \ref{cond:unifbd} is satisfied.

\vspace{2em}

\bibliographystyle{apalike}
\bibliography{Citation.bib}
\end{document}